\documentclass[11pt,a4paper]{article}
\usepackage{ifthen}
\usepackage{amsmath}
\usepackage{amssymb}
\usepackage{latexsym}
\usepackage{bbm}
\usepackage{graphicx}
\usepackage{pdfsync}
\usepackage{verbatim}
\usepackage{fancyhdr}
\usepackage{url}               
\usepackage[english]{babel}
\usepackage[applemac]{inputenc}

\usepackage[pdftex,%
a4paper=true,%
bookmarks=true,%
bookmarksnumbered=true,%
colorlinks=true,%
urlcolor=blue,%
pdfstartview=FitH%
]{hyperref}

\newtheorem{theorem}{Theorem}[section]
\newtheorem{lemma}[theorem]{Lemma}

\newtheorem{proposition}[theorem]{Proposition}
\newtheorem{corollary}[theorem]{Corollary}

\newtheorem{remark}[theorem]{\it \bf{Remark}\/}

\numberwithin{equation}{section}
\numberwithin{figure}{section}

\newcommand{\bp}{{\it Proof. }}
\newcommand{\ep}{\hfill $\square$\\}

\newcommand{\bpl}{{\it Proof of Lemma }} 
\newcommand{\epl}{\hfill $\square$\\}

\newcommand{\bpp}{{\it Proof of Proposition }} 
\newcommand{\epp}{\hfill $\square$\\}

\newcommand{\Imm}{\mbox{Im}}
\renewcommand{\Re}{\mbox{Re}}
\newcommand{\be}{\begin{equation}}
\newcommand{\ee}{\end{equation}}
\newcommand{\bea}{\begin{eqnarray}}
\newcommand{\eea}{\end{eqnarray}}
\newcommand{\bee}{\begin{eqnarray*}}
\newcommand{\eee}{\end{eqnarray*}}

\def\pa{\partial}

\def\CC{\mathbb{C}}
\def\NN{\mathbb{N}}
\def\RR{\mathbb{R}}

\def\dps{\displaystyle}
\def\ni{\noindent}

\def\eps{\vare}

\def\eps{\varepsilon}

\catcode`@=11
\def\supess{\mathop{\operator@font Sup\,ess}}
\catcode`@=12
\def\CC{\mathbb{C}}

\def\NN{\mathbb{N}}

\def\RR{\mathbb{R}}
\def\CC{\mathbb{C}}

\def\KK{\mathbb{K}}

\def\ds{\displaystyle}

\def\ni{\noindent}

\def\bar#1{{\overline #1}}

\def\R2+{\RR ^2_+}

\def\pa{\partial}

\def\lim{\mathop{\rm lim}}

\def\log{{\rm log}}

\def\pa{\partial}

\def\pa{\partial}

\begin{document}

\renewcommand{\refname}{References}
\bibliographystyle{alpha}

\pagestyle{fancy}
\fancyhead[L]{ }
\fancyhead[R]{}
\fancyfoot[C]{}
\fancyfoot[L]{ }
\fancyfoot[R]{}
\renewcommand{\headrulewidth}{0pt} 
\renewcommand{\footrulewidth}{0pt}

\newcommand{\montitre}{Construction of an eigen-solution for the Fokker-Planck operator with heavy tail equilibrium: an \`a la Koch method in dimension 1}

\newcommand{\auteur}{\textsc{Dahmane Dechicha$^\ast$,  Marjolaine Puel$^\ddagger$ }}
\newcommand{\affiliation}{$^\ast$Laboratoire J.-A. Dieudonn\'e.  Universit\'e C\^{o}te d'Azur. UMR 7351 \\ Parc Valrose, 06108 Nice Cedex 02,   France\\
$^\ddagger$Laboratoire de recherche AGM. CY Cergy Paris Universit\'e. UMR CNRS 8088 \\

2 Avenue Adolphe Chauvin 95302 Cergy-Pontoise Cedex, France \\
\url{dechicha@unice.fr}; \url{ mpuel@cyu.fr
} }

 \begin{center}
{\bf  {\LARGE \montitre}}\\ \bigskip \bigskip
 {\large\auteur}\\ \bigskip \smallskip
 \affiliation \\ \bigskip
\today
 \end{center}
 \begin{abstract}
This paper is devoted to the construction of an \emph{eigen-solution} for the Fokker-Planck operator with heavy tail equilibrium. We propose an \textit{alternative} method in dimension 1, which will be generalizable in higher dimension. The later method is inspired by the work of H. Koch on non-linear KdV equation \cite{Koch}. As a consequence of this construction, we recover the result of G. Lebeau and M. Puel \cite{LebPu} on the fractional diffusion limit for the Fokker-Planck equation. 
\end{abstract}

\tableofcontents
\newpage

 \pagestyle{fancy}
\fancyhead[R]{\thepage}
\fancyfoot[C]{}
\fancyfoot[L]{}
\fancyfoot[R]{}
\renewcommand{\headrulewidth}{0.2pt} 
\renewcommand{\footrulewidth}{0pt} 

\section{Introduction}
\subsection{Setting of the problem}
Our starting point is the kinetic Fokker-Planck (FP) equation, which describes in a deterministic way the \emph{Brownian motion} of a set of particles. It is given by the following form
\begin{equation}\label{FP}
\left\{\begin{array}{l}
 \partial_t f + v\cdot \nabla_x f = Q(f), \; \; \; t \geqslant 0 , \ x \in \RR^d , \ v \in \RR^d ,\\
\\
    f(0,x,v) = f_0(x,v), \; \; \quad x \in \RR^d , \ v \in \RR^d ,
\end{array}\right.
\end{equation}
where the collisional Fokker-Planck operator $Q$ is given by
\begin{equation}\label{defQ}
Q(f)=\nabla_v\cdot\bigg(F\nabla_v\bigg(\frac{f}{F}\bigg)\bigg) ,
\end{equation}
and $F$ is the equilibrium of $Q$, a fixed function which depends only on $v$ and satisfying  
$$ Q(F)=0 \quad \mbox{ and } \quad \int_{\RR^d} F(v)\ \mathrm{d}v = 1 . $$ 
For non-negative initial data $f_0$, the unknown $f(t,x,v) \geqslant 0$ can be interpreted as the density of particles occupying at time $t \geqslant 0$, the position $x \in \RR^d$ with velocity $v \in \RR^d$. \\

The \emph{classical} or \emph{fractional diffusion limits} for kinetic equations have been studied in a series of papers in recent years. The principal motivation behind this study is to derive simpler models from kinetic equations with collision operators. Indeed, in the diffusion approximation, the velocity variable is involved only in the equilibrium. This approximation is a kind of hydrodynamic limit. When the interactions between particles are the dominant phenomena and when the time of observation is very large, it reflects the fact that we are between the mesoscopic and the macroscopic scale. More precisely, we introduce a small parameter $\eps \ll 1$, the mean free path and we proceed to rescaling the distribution function $f(t,x,v)$ in time and space
$$ t=\frac{t'}{\theta(\eps)} \quad \mbox{ and } \quad x=\frac{x'}{\eps} \quad \mbox{ with } \quad \theta(\eps) \underset{\eps \to 0}{\longrightarrow } 0 , $$
which leads to the following rescaled equation (without primes)
\begin{equation}\label{fp-theta}
\left\{\begin{array}{l}
\theta(\varepsilon) \partial_t f^\varepsilon + \varepsilon v\cdot \nabla_x f^\varepsilon = Q(f^\varepsilon), \; \; \;  t \geqslant 0 , x \in \RR^d , v \in \RR^d , \\
\\
    f^\varepsilon(0,x,v) = f_0(x,v), \; \; \; \; x \in \RR^d , v \in \RR^d .
\end{array}\right.
\end{equation}
Note that taking initial conditions independent of $\varepsilon $ means that we take well prepared initial conditions in the non rescaled variable. Now, in order to study the behaviour of the solution $f^\eps$ as $\eps \to 0$, we are interested in the following question: what is the good time scale $\theta(\eps)$ such that the limit of the solution $f^\eps$ is not trivial?  Formally, passing to the limit when $\eps \to 0$ in the equation \eqref{fp-theta}, we obtain that the limit $f^0$ is in the kernel of $Q$ which is spanned by the equilibrium $ F $, which means that $f^0=\rho(t,x)F(v)$. Thus, what is the equation satisfied by $\rho$ and how do we justify this passage \emph{rigorously}? \\

The limit of $f^\eps$ when $\eps$ goes to $0$ may depend on the nature of the equilibrium considered. For Gaussian equilibria, the answers have been known for a long time and it is the \emph{classical diffusion} that we obtain for a usual scaling $\theta(\eps)=\eps^2$. For slowly decreasing equilibria, or so-called \emph{heavy-tailed equilibria} of the form $F(v) \sim \langle v \rangle^{-\beta}$, it is more complicated, and the answers to this question has been the interest of many papers in the last few years, with different methods and for different collision operators. This was initiated by Mischler, Mouhot and Mellet \cite{MMM} on the linear Boltzmann equation, for an ad hoc cross section, which does not depend on the spatial variable, where they obtained classical and fractional diffusion, according to the powers of $\langle v \rangle$ which appear in $F$. Their method is based on \emph{Fourier-Laplace transformation}, with close links to earlier work by Milton, Komorowski and Olla \cite{MKO} on \emph{Markov chains}. This was continued by A. Mellet \cite{M}, still for the same equation, more general since it also applies to cross sections that depend on the position variable. See also the work of Ben Abdallah, Mellet and the second author \cite{BMP AD, BMP FD}, where they used a \emph{Hilbert expansions approach} \cite{BMP FD} and obtained a strong convergence of the solution $f^\eps$ to $\rho F$ for initial data $f_0$ in $H^4(\RR^{2d})$. \\

The diffusion limit for the FP equation has been studed after that of the linear Boltzmann one,  and seems more complicated.  Indeed,  there is no \emph{spectral gap}, thus makes the limiting process more difficult. In addition, for this equation, all the terms of the operator participate in the limit, i.e. the collision and advection parts.  The limit was studed in \cite{NP} with $\theta(\eps)=\eps^2$ and $\beta>d+4$ using the \emph{moment method}. Also in \cite{CNP} for the critical case $\beta=d+4$ with $\theta(\eps)=\eps^2|\log \eps|$ using \emph{a probabilistic approach}. In both cases, the authors obtained classical diffusion.  Concerning the case $d<\beta<d+4$, it has been solved recently, by G. Lebeau and the second author in dimension 1 \cite{LebPu}, with a \emph{spectral approach}, then by N. Fournier and C. Tardif in dimension 1 \cite{FT d1} and then in dimension $d \geqslant 2$ \cite{FT} by \emph{a probabilistic approach}. Depending on the range of the exponents, different regimes corresponding to Brownian processes, stable processes or integrated symmetric Bessel processes are obtained and described in this last paper. Recently, E. Bouin and C. Mouhot \cite{BM} have constructed what they call a `\emph{fluid mode}' using a method that combines energy estimates and a quantitative spectral method. Their method was valid to scattering models, Fokker-Planck (by introducing a \emph{weighted} density, $\ \tilde \rho := \int_{\RR^d} \frac{f}{\langle v\rangle^2} \mathrm{d}v$) and L\'evy-Foker-Planck operators. We refer also to the paper by E. Bouin, J. Dolbeau and L. Lafleche \cite{BDL}, where the authors have developed an $L^2$-hypocoercivity approach and established an optimal decay rate, determined by a fractional Nash type inequality, compatible with the fractional diffusion limit.  \\

In \cite{LebPu}, the authors proposed to take the eigenfunction of the Fokker-Planck operator (with the advection part) as an adequate test function, i.e. take the solution of the spectral problem
 $$ (Q+\mathrm{i}\eps \xi v)M_{\eps,\mu}(v) = \mu M_{\eps,\mu}(v) , \quad v \in \RR ,$$
written in Fourier variable, $\xi$, for the position $x$. Their construction consists in reconnecting two branches constructed as follows: They first constructed for each $\mu, \eps$ (fixed) a branch in the half space $\RR_+$, by introducing an approximate equation for large velocities. Then, by symmetry of the equation, they obtained a second branch in the other half space $\RR_-$, and in order to reconnect the two functions to obtain a $C^1(\RR)$ solution, the reconnection of the derivatives implies a relation $\mu(\eps)$.  This method of reconnection seems very complicated to adapt in higher dimension, because it means that the implicit function theorem; used to study constraint; must be applied to the whole derivative operator.  \\

The purpose of this paper is to propose an ``\emph{alternative}'' method, inspired by the work of Herbert Koch on the non-linear KdV equation \cite{Koch}, where we solve the spectral problem associated to the Fokker-Planck operator, taking into account the advection part, working on the whole space $\RR$ and avoiding the reconnection problems. This method is generalizable in higher dimension and is probably interesting for generalized potential or non linear problems.  \\
 
In our proof, as we will explain later, a splitting of the Fokker-Planck operator is involved, which recalls the \emph{enlargement theory} for nonlinear Boltzmann operator when there are spectral gap issues. This theory was developed by Gualdani, Miscler and Mouhot in \cite{GMM} whose key idea was based on the decomposition of the operator into two parts, a \emph{dissipative} part plus a \emph{regularizing} part. See also P. Gervais (\cite{Gervais} and references therein)  for a spectral study of the linearized Boltzmann operator in $L^2$ spaces with polynomial and Gaussian weights.

\subsection{Setting of the result}
In the present work, we consider for any $\beta > 1$, heavy tail equilibria
$$ F(v) \sim \frac{1}{(1+|v|^2)^\frac{\beta}{2}} .$$
Before stating our main result, let us give some notations that we will use along this paper.\\

\noindent \textbf{Notations.}
As in \cite{LebPu}, in order to simplify the computation and work with a self-adjoint operator in $L^2$, we proceed to a change of unknown by writing
$$ f = F^\frac{1}{2}g = C_\beta Mg $$
with  $$ M=C_\beta^{-1} F^\frac{1}{2}=\frac{1}{(1+|v|^2)^\frac{\gamma}{2}} , $$ 
since we impose $\gamma :=\frac{\beta}{2}>\frac{1}{2}$, $F \in L^1(\RR)$ then, $M \in L^2(\RR)$ and $C_\beta$ is chosen such that $$ \int_{\RR} F dv = 1 .$$
The equation \eqref{fp-theta} becomes
$$ \theta(\eps) \partial_t g^\eps + \eps v\cdot \nabla_x g^\eps =\frac{1}{M}\pa_v \bigg(M^2 \pa_v\bigg(\frac{g^\eps}{M}\bigg)\bigg) = \pa_v^2 g^\eps - W(v)g^\eps ,
$$
with 
$$ W(v)= \frac{\pa_v^2 M}{M} = \frac{\gamma(\gamma+1)|v|^2 -\gamma}{(1+|v|^2)^2}.
$$
We see the equation as 
$$
\theta(\eps) \partial_t g^\eps = - \mathcal{L}_\eps g^\eps ,
$$
where
$$ \mathcal{L}_\eps := -\pa_v^2 + W(v) + \eps v \cdot \nabla_x = -(Q-\eps v \cdot \nabla_x)$$
and $$ Q =  -\Delta_v + W(v) .$$
We operate a Fourier transform in $x$ and since the operator $Q$ has coefficient that do not depend on $x$, we get:
\begin{equation}\label{rescaled}
\theta(\eps) \partial_t \hat g^\eps = - \mathcal{L}_\eta \hat g^\eps ,
\end{equation}
where
$$ \mathcal{L}_\eta := -\pa_v^2 + W(v) + \mathrm{i} \eta v $$
and $$ \eta = \eps \xi ,$$ 
with $\xi$ being the space Fourier variable. \\

\noindent The operator $\mathcal{L}_\eta$ is an unbounded self-adjoint operator acting on $L^2$. Its domain is given by
$$D(\mathcal{L}_\eta) = \big\{  g \in L^2(\mathbb{R}) \ ; \ \pa_v^2 g \in L^2(\mathbb{R}), v g \in L^2(\mathbb{R})\big\} . $$

The aim of this paper is to prove, with a geometry independent method, the following main theorem.
\subsubsection*{Main Theorem}
\begin{theorem}[Eigen-solution for the Fokker-Planck operator]\label{main}
 Assume that $\beta \in ]1,5[\setminus\{2\}.$ Let $\eta_0>0$ and $\lambda_0>0$ small enough.  Then, for all $\eta\in [0,\eta_0]$,  there exists a unique eigen-couple $\big(\mu(\eta),M_\eta\big)$ in $\{\mu\in \mathbb C, |\mu|\leqslant \eta^{\frac{2}{3}}\lambda_0\}\times L^2(\RR,\CC)$, solution to the spectral problem 
\begin{equation}\label{M_mu,eta}
\mathcal{L}_\eta(M_{\mu,\eta})(v)=\big[-\partial_v^2 +W(v) + \mathrm{i} \eta v \big] M_{\mu,\eta} (v)= \mu M_{\mu,\eta}(v) , \  v\in \RR .
\end{equation}
Moreover,  one has 
\begin{enumerate}
\item The following convergence in $L^2(\RR,\CC)$,
\begin{equation}\label{item 1 thm}
 \| M_\eta - M \|_{L^2} \underset{\eta \to 0}{\longrightarrow} 0 .
\end{equation} 
\item The relationship between the eigenvalue $\mu(\eta)$, the scale of the time variable $\theta(\eps)$ and the coefficient $\kappa$ is given by the following expansion:
\begin{equation}\label{glres}
\mu(\eta) = \kappa \eta^{\frac{\beta+1}{3}}\big(1+O(\eta^{\frac{\beta+1}{3}})\big) ,
\end{equation}
where $\kappa$ is a positive constant given by 
\begin{equation}\label{kappa thm}
\kappa= -2C_\beta^2 \int_0^\infty s^{1-\gamma} \Imm H_0(s) \mathrm{d}s ,
\end{equation}
and where $H_0$ is the unique solution to the equation 
\begin{equation}
\big[-\pa_s^2+\frac{\gamma(\gamma+1)}{s^2} + \mathrm{i} s \big]H_0(s) = 0 , \quad \forall s \in \RR^* ,
\end{equation}
satisfying 
\begin{equation}
\int_{\{|s|\geqslant 1 \}} |H_0(s)|^2 \mathrm{d}s < \infty \quad \mbox{ and } \quad H_0(s) \underset{0}{\sim} |s|^{-\gamma} .
\end{equation}
\end{enumerate}
For $\eta\in [-\eta_0,0]$, by complex conjugation on the equation, we get 
$$\mu(\eta)=\bar\mu(-\eta)= \kappa |\eta|^{\frac{\beta+1}{3}}(1+O(|\eta|^{\frac{\beta+1}{3}})).$$
\end{theorem}
\begin{remark}The hypothesis $\beta\neq 2$ is  technical. It avoids to introduce logarithmic terms
in the expression of $\mu(\eta)$. 
\end{remark}
\subsubsection*{Idea of the proof}
The proof of our main result is done in two main steps, both based on the Implicit Function Theorem.
First, we consider what we call a \emph{penalized equation}.  We introduce an additional term that kills the $M$ direction in the kernel of the linear operator computed at $\eta=0$. That gives the following equation
\begin{equation}\label{eq penalisee1}
\left\{ \begin{array}{l}
\big[-\partial_v^2 +W(v) + \mathrm{i} \eta v \big] M_{\mu,\eta} (v)= \mu M_{\mu,\eta}(v) - \langle M_{\mu,\eta}-M,\Phi\rangle \Phi , \  v\in \mathbb{R} ,\\
\\
M_{\mu,\eta} \in L^2(\mathbb{R}).
\end{array}\right.
\end{equation}
where $\Phi$ is a function that we will determine later.  This penalized equation has a solution for any $\lambda$ and $\eta$ on the whole space and it allows us to avoid the problem of reconnection and to work directly on the whole space $\RR$. This is one of the key points of this method and it iallows to generalize this construction in any dimension. \\

The objective of the first step is to show the existence of a unique solution for equation \eqref{eq penalisee1} for any  $\eta$ and $\mu$. Indeed, as we said above, we will decompose the operator $``-\partial_v^2 +W(v) + \mathrm{i} \eta v -  \mu"$ into two parts. The first part is chosen such that it admits ``a right inverse" that is continuous as a linear operator between two suitable functional spaces, continuous with respect to the parameters $\eta$ and $\mu$ and compact at $\eta=\mu=0$.  The second part of the operator is left in the right-hand side of the equation, i.e. is considered as a source term. Find a solution to \eqref{eq penalisee1} remains to a fixed point process. \\

In the second step, to ensure that the additional term vanishes, we have to chose $\mu(\eta)$ obtained via the Implicit Function Theorem around the point $(\mu,\eta)=(0,0)$.  

\subsection{Relation to the fractional diffusion problem}
In this subsection, we will explain how one can recover the fractional diffusion limit for the Fokker-Planck equation \eqref{fp-theta}, and how the eigenvalue is related to the diffusion coefficient in our main result. See Section 3 of \cite{LebPu} for more details. \\

\noindent \textbf{Heuristic on the computation of the eigenvalue.} With a formal calculation, we will present how the time scaling $\theta(\eps)$ is chosen and how it appears in the spectral problem. Assume that the couple $(\mu(\eta),M_{\mu,\eta})$ is solution to the problem
$$
\mathcal{L}_\eta (M_{\mu,\eta}) = [-\pa_v^2+W(v)+\mathrm{i} \eta v]M_{\mu,\eta} = \mu(\eta) M_{\mu,\eta} , \quad v \in \RR.
$$
Then, integrating this equation against $M$ and using the fact that $[-\pa_v^2+W(v)]M=0$, we get
$$
\mathrm{i} \eta \int_\RR v M_{\mu,\eta} M \mathrm{d}v = \mu(\eta) \int_\RR M_{\mu,\eta} M \mathrm{d}v .$$
Therefore,
$$
 \mu(\eta) = - \mathrm{i} \eta \int_\RR v M_{\mu,\eta} M \mathrm{d}v \bigg(\int_\RR M_{\mu,\eta} M \mathrm{d}v\bigg)^{-1} .$$
 If $M_{\mu,\eta} \to M $ when $\eta \to 0$ then, we get $\big(\int_\RR M_{\mu,\eta} M \mathrm{d}v\big)^{-1} \to \|M\|_2^{-2}=C_\beta^2$. Formally, by a Hilbert expansion 
 $$
 M_{\mu,\eta} =M+\eta N+o(\eta), \mbox{ where } [-\pa_v^2+W(v)]N=+\mathrm{i} v M,
 $$
 we get $$|M_{\mu,\eta} (v) - M(v) | \lesssim  \eta \langle v \rangle^{3-\gamma}  \quad \mbox{ for } \ |v|\leqslant \eta^{-\frac{1}{3}}s_0 ,$$ 
and by rescaling the integral on large velocities by $v=\eta^{-\frac{1}{3}} s$ we finally obtain:
 $$  \mu(\eta) \sim \eta^2 \int_{|v|\leqslant \eta^{-\frac{1}{3}}s_0} v \langle v \rangle^{4-2\gamma} \mathrm{d}v + \eta^{\frac{2\gamma+1}{3}} \int_{s \geqslant s_0} s^{1-\gamma} \Imm H_0(s) \mathrm{d}s ,$$
where $H_0$ is the limit of $M_{\mu,\eta}$ rescaled,  solution to the limit ``rescaled equation''
$$ \big[ -\pa_s^2 + \frac{\gamma(\gamma+1)}{s^2} + \mathrm{i}  s\big]H_0(s) = 0 , \quad \forall s \in \RR^* .$$
Recall that $\eta=\eps \xi$. Then, $\mu(\eta)=\mu(\eps|\xi|)$. Thus, if $2\gamma>5$ then we find classical diffusion with the usual scaling $\theta(\eps)=\eps^2$, and it would be the small velocities which give the diffusion coefficient. While if $2\gamma<5$ then we get fractional diffusion with a power of $\frac{2\gamma+1}{6}$ for the Laplacian, the scaling is given by $\theta(\eps)=\eps^{\frac{2\gamma+1}{3}}$ and the diffusion coefficient is determined by the integral for large velocities in this case.  \\

\noindent {\bf Heuristic of the diffusion approximation}

\begin{remark}
 Note that we will work with the Fourier transform of $\rho$ and we will establish that 
$ \hat\rho(t,\xi)=\int e^{-\mathrm{i}x\xi} \rho(t,x)\mathrm{d}x$ satisfies
\begin{equation}\label{diff}
\partial_t\hat \rho+\kappa|\xi|^{\alpha}\hat\rho=0.
\end{equation}
\end{remark}

\noindent Let $\xi \in \mathbb R$ and let $M_\eta$ be the unique solution in $L^2$ of the equation $\mathcal{L}_\eta(M_\eta)=\mu(\eta)M_\eta $ given in Theorem \ref{main}. One has
$$
\begin{array}{rcl}
\dps \frac{\partial}{\partial t} \int  \hat g^\eps(t,\xi,v) M_\eta \mathrm{d}v&=&\dps  \int \partial_t\hat g^\eps M_\eta \mathrm{d}v
= -\eps^{-\alpha}\int \mathcal{L}_\eps(\hat g^\eps) M_\eta \mathrm{d}v\\
&=&
\dps -\eps^{-\alpha}\int \hat g^\eps    \mathcal{L}_\eps(M_\eta) \mathrm{d}v = -\eps^{-\alpha}\mu(\eta)\int \hat g^\eps(t,\xi,v)  M_\eta \mathrm{d}v  .  
\end{array}
$$
Recall that $\eta = \eps \xi$ and by Theorem \ref{main}, $\eps^{-\alpha}\mu(\eta) \to \kappa |\xi|^\alpha$ and $M_\eta \to M$ when $\eps \to 0$.
Therefore, passing to the limit in the previous equality formally, we obtain the equation
$$ \partial_t \hat \rho = |\xi|^\alpha\kappa \hat\rho.$$
Which means that $\hat \rho$ satisfies \eqref{diff}. Thus, the solution $f^\eps$ of \eqref{fp-theta} converges; in some sens; towards $\rho(t,x)F(v)$, where $\rho$ is the solution of the fractional diffusion equation
\begin{equation}
\partial_t\rho +\kappa (-\Delta)^{\frac{\alpha}{2}}\rho =0,\quad \rho(0,x)=\int f_0 \mathrm{d}v  .
 \end{equation}
Returning to the space variable, the fractional Laplacian $(-\Delta)^{\frac{\alpha}{2}}$ is a non-local operator which can be defined by
$$ (-\Delta)^{\frac{\alpha}{2}}\rho(x) := c_{\alpha} \ \mathrm{P.V.} \int_{\RR} \frac{\rho(x)-\rho(y)}{|x-y|^{\alpha+1}} \mathrm{d}y , $$
which also can be seen as an $\alpha$-stable L\'evy process, and thus can be interpreted as a random trajectory, generalising the concept of Brownian motion, which may contain jump discontinuities.

\subsection{Notations and definition of the considered operators}
Let $\lambda$ and $\eta$  be fixed, where $\lambda \in \CC$ being defined by $\lambda := \mu \eta^{-\frac{2}{3}}$. Let us denote by $L_{\lambda,\eta}$ the operator
$$ L_{\lambda,\eta} := -\partial_v^2 + \tilde W(v) + \mathrm{i} \eta v -  \lambda \eta^\frac{2}{3} ,$$
where $$\tilde W(v)=\frac{\gamma(\gamma+1)}{1+v^2},$$
and let denote by $V := \tilde W - W$.  We have
$$ V(v)=\frac{\gamma(\gamma+2)}{(1+v^2)^2}.$$
We will rewrite the equation \eqref{eq penalisee1} as follows:
\begin{equation}\label{eq penalisee2}
\left\{ \begin{array}{l}
L_{\lambda,\eta} M_{\lambda,\eta} (v) =  V(v) M_{\lambda,\eta}(v) - \langle M_{\lambda,\eta}-M,\Phi\rangle \Phi(v) , \ v\in \mathbb{R} ,\\
\\
M_{\lambda,\eta} \in L^2(\mathbb{R}).
\end{array}\right.
\end{equation}
The two equations \eqref{eq penalisee1} and \eqref{eq penalisee2} are equivalent.
\begin{remark}\label{symetrie de la sol}
\item \begin{enumerate}
\item Since $L_{\lambda,0}$ does not depend on $\lambda$, let's denote it by $L_0$, $L_0:=L_{\lambda,0}$.
\item If $ \ \bar \Phi(-v)=\Phi(v)$ and $M_{\lambda,\eta}(v)$ satisfies the equation \eqref{eq penalisee2}, then $\overline{M}_{\bar \lambda,\eta}(-v)$ is also solution to \eqref{eq penalisee2}, since the potential $W$ is symmetric for a symmetric equilibrium $M$. Note that this is where the symmetry of the equilibrium $M$ is used and therefore this is a ``non-drift condition".
\item In the first step of our proof, the continuity of $L_{\lambda,\eta}^{-1}$ and the compactness of $L_0^{-1}$ are among the most important and difficult points to prove. The splitting of the potential $W$ into $V$ and $\tilde W$ is the key idea in the proof of these last two points.
\end{enumerate}
\end{remark}

\noindent We will construct a right inverse to the operator $L_{\lambda,\eta}$ defined above. We are going to study the operator $L_0$ first, then we construct solutions for the equation $L_{\lambda,\eta}(\psi)=f$, with any source term $f$, from the solutions of $L_0(\psi)=f$, after noticing that $L_{\lambda,\eta}$ is in fact a small perturbation of $L_0$ for $\lambda$ and $\eta$ small enough.

\section{Green function for the leading part of the limiting operator }
We consider the equation:
\begin{equation}\label{L_0(psi)=f}
 L_0(\psi)=\big[-\partial^2_v+\tilde W(v)\big]\psi=f, \quad v \in \RR
\end{equation}
We are going to construct solutions for this equation by the method of the variation of the constant, which allows us to define a right inverse to $L_0$. Let's start by looking at the homogeneous equation: 
\begin{equation}\label{L_0(psi)=0}
L_0(\psi)=\big[-\partial^2_v+\tilde W(v)\big]\psi=0 , \quad v \in \RR.
\end{equation}

\subsection{Construction of an intermediate solution to the homogeneous equation}
In order to construct the basis of solution to the homogeneous equation, we need the following two lemmas:
\begin{lemma}[Behaviour of the solutions of \eqref{L_0(psi)=0}]\label{comportement psi}
If $\tilde \psi$ is a solution for the equation \eqref{L_0(psi)=0} then, either $\tilde \psi(v) \underset{+\infty}{\sim} |v|^{-\gamma}$ or $\tilde \psi(v) \underset{+\infty}{\sim} |v|^{\gamma+1}$. Thus, for a positive constant $v_0$ large enough, there exists a unique function $R_1 \in L^\infty([v_0,+\infty[)$, with $R_1=O(v^{-2})$ such that 
$$ \tilde \psi= |v|^{-\gamma}(1+R_1) \quad \mbox{ or } \quad \tilde \psi=|v|^{\gamma+1}(1+R_1) . $$
Similarly, either $\psi(v) \underset{-\infty}{\sim} |v|^{-\gamma}$ or $\psi(v) \underset{-\infty}{\sim} |v|^{\gamma+1}$. Then, for a positive constant $v_0$ large enough, there exists a unique function $R_2 \in L^\infty(]-\infty,-v_0])$, with $R_2=O(|v|^{-2})$) such that
$$\tilde \psi = |v|^{-\gamma}(1+R_2) \quad \mbox{ or } \quad \tilde \psi = |v|^{\gamma+1}(1+R_2).$$
\end{lemma}

\noindent \begin{lemma}[The existence of an intermediate solution for \eqref{L_0(psi)=0}]\label{existence de psi} 
\item There is a unique positive smooth function $\psi$, solution to the Cauchy problem:
\begin{equation}\label{pb de Cauchy psi}
\left\{\begin{array}{l}
 L_0(\psi)=0, \\
   \psi(0)=1, \\
   \psi'(0)=0.
\end{array}\right.
\end{equation}
Moreover, 
\begin{enumerate}
\item For all $v$ in $\RR$: $\psi(-v)=\psi(v)>0$. 
\item There exists a positive constant $v_0$ large enough, there exists a unique function $R$ belongs to $L^\infty([v_0,+\infty[)$, with $R(v)=O(|v|^{-2})$ such that 
\begin{equation}
\psi(v)=c|v|^{\gamma+1}\big(1+R(v)\big), \quad \forall |v|\geqslant v_0 ,
\end{equation}
where $c$ is a given positive constant.
\end{enumerate} 
\end{lemma}
\bpl \ref{comportement psi}. Since $\tilde W(v)\underset{|v|\sim \infty}{\sim}\frac{\gamma(\gamma+1)}{v^2}$ for $v\gg 1$, we consider the approximate equation 
\begin{equation}\label{tilde L_0}
\tilde L_0(\psi):=\big[-\partial_v^2+\frac{\gamma(\gamma+1)}{v^2}\big]\psi=0 .
\end{equation}
The two functions $\psi_1^*=|v|^{-\gamma}$ and $\psi_2^*=|v|^{\gamma+1}$ are solutions to \eqref{tilde L_0} in $\RR^*$, and form a basis of solutions on $\RR^*_+$ and $\RR^*_-$ respectively. Now, since we know the asymptotic profile of the solutions to \eqref{L_0(psi)=0}, on $\{|v|\geqslant v_0\}$, we will look for solutions via the following change of unknown 
$$\tilde \psi_1= \psi_1^*(1+R_1) \mbox{ on } [v_0,+\infty) \ \mbox{ and } \ \tilde \psi_2= \psi_2^*(1+R_2) \mbox{ on } (-\infty,-v_0] .$$
We show the existence and uniqueness of $R_i$ for $i = 1, 2$ by a fixed point argument performed on $R_i$ which satisfy the equation
$$2(\psi_i^*)'R_i'+\psi_i^* R_i''=[\tilde W - \frac{\gamma(\gamma+1)}{v^2}]\psi_i^*(1+R_i) ,$$
which can be written (after multiplication by $\psi_i^*$) as follows
$$ 
\big([\psi_i^*]^2 R_i'\big)'=[\tilde W - \frac{\gamma(\gamma+1)}{v^2}][\psi_i^*]^2(1+R_i) 
$$
and this leads to the implicit equations
$$
\left\{\begin{array}{l}
\ds R_1(v)= \int_v^\infty\int_v^w \frac{[\psi_j^*]^2(w)}{[\psi_j^*]^2(u)}\mathrm{d}u\big[\tilde W(w) - \frac{\gamma(\gamma+1)}{w^2}\big](1+R_1(w))\mathrm{d}w , \ v\geqslant v_0 , \ j=1,2 ,
\\
\\
 \ds R_2(v)=\int_{-\infty}^v\int_w^v \frac{[\psi_j^*]^2(w)}{[\psi_j^*]^2(u)}\mathrm{d}u\big[\tilde W(w) - \frac{\gamma(\gamma+1)}{w^2}\big](1+R_2(w))\mathrm{d}w , \ v\leqslant -v_0 , \ j=1,2 .
\end{array}\right.
$$
Note that the two previous equations are of the form: $(Id-\mathbb{K}_i)(R_i)=\mathbb{K}_i(1)$ for $i=1, 2$ with
$$
\left\{\begin{array}{l}
\ds \mathbb{K}_1(g)= \int_v^\infty\int_v^w \frac{[\psi_j^*]^2(w)}{[\psi_j^*]^2(u)}\mathrm{d}u\big[\tilde W(w) - \frac{\gamma(\gamma+1)}{w^2}\big]g(w)\mathrm{d}w , \ v\geqslant v_0 , \  j=1,2 ,
\\
\\
 \ds \mathbb{K}_2(g)=\int_{-\infty}^v\int_w^v \frac{[\psi_j^*]^2(w)}{[\psi_j^*]^2(u)}\mathrm{d}u\big[\tilde W(w) - \frac{\gamma(\gamma+1)}{w^2}\big]g(w)\mathrm{d}w , \ v\leqslant -v_0 ,  \ j=1,2 .
\end{array}\right.
$$
We have for $v\geqslant v_0 > 0$, \\
$\bullet$ for $j=1$:
\begin{align*}
\big|v^{n+2}\mathbb{K}_1(g)(v)\big|&\leqslant v^{n+2}\int_v^\infty\int_v^w \frac{u^{2\gamma}}{w^{2\gamma}}\mathrm{d}u\frac{\gamma(\gamma+1)}{w^{n+4}}w^n|g(w)|\mathrm{d}w \\
&\leqslant \frac{\gamma(\gamma+1)}{2\gamma+1}v^{n+2}\int_v^\infty w^{-3-n}\mathrm{d}w \|v^n g\|_{L^\infty(v_0,\infty)} \\
&\leqslant \frac{\gamma(\gamma+1)}{2\gamma+1}\frac{1}{n+2} \|v^n g\|_{L^\infty(v_0,\infty)} .
\end{align*}
$\bullet$ for $j=2$, we take $n>2\gamma$ and we write:
\begin{align*}
\big|v^{n+2}\mathbb{K}_1(g)(v)\big|&\leqslant v^{n+2}\int_v^\infty\int_v^w \frac{w^{2\gamma+2}}{u^{2\gamma+2}}\mathrm{d}u\frac{\gamma(\gamma+1)}{w^{n+4}}w^n|g(w)|\mathrm{d}w \\
&\leqslant \frac{\gamma(\gamma+1)}{2\gamma+1}v^{n-2\gamma-1}\int_v^\infty w^{2\gamma-n-2}\mathrm{d}w \|v^n g\|_{L^\infty(v_0,\infty)} \\
&\leqslant \frac{\gamma(\gamma+1)}{2\gamma+1}\frac{1}{n+2-2\gamma} \|v^n g\|_{L^\infty(v_0,\infty)}  .
\end{align*}
Finally,  in both cases, $\mathbb{K}_1$ is bounded in $L^\infty(v_0,\infty)$ with
$$ \|\mathbb{K}_1\|_{\mathcal{L}(L^\infty(v_0,\infty))}\leqslant \frac{\gamma(\gamma+1)}{2(2\gamma+1)} v_0^{-2}\leqslant\frac{1}{2} \quad \mbox{ for } v_0\geqslant\sqrt{\frac{\gamma(\gamma+1)}{2(2\gamma+1)}}.$$
Hence the existence and uniqueness of $R_1$ in $L^\infty(v_0,\infty)$ with the asymptotic $R_1=O(v^{-2})$ since $\mathbb{K}_1(1)=O(v^{-2})$. \\
By doing exactly the same computations for $R_2$, as in $R_1$ (after change of variable $v=-v'$), we get
$$
\big||v|^{n}\mathbb{K}_2(g)(v)\big| \leqslant \frac{\gamma(\gamma+1)}{2(2\gamma+1)} v_0^{-2}\big\||v|^n g \big\|_{L^\infty(-\infty,-v_0)}
, \ \mbox{ for all } j =1,2  \ \mbox{ if } n> 2\gamma . $$
The rest of the proof is done in the same way as for $R_1$.
\epl
\begin{remark}
Note that $0$ is the unique solution to \eqref{L_0(psi)=0} in $L^2(\RR)$. It can be seen by multiplying the equation by $\psi$ and integrating by parts over $\RR$. Therefore, we cannot have a solution like $|v|^{-\gamma}$ in $+\infty$ and $-\infty$ at the same time. 
\end{remark}

\noindent \bpl \ref{existence de psi}. \\
\textbf{Existence:} The existence of a unique global solution $\psi$ for the Cauchy problem \eqref{pb de Cauchy psi} follows from the Cauchy-Lipschitz theorem.
Concerning regularity, we have $\partial_v^2 \psi = \tilde W \psi \in C^0(\RR)$, then $\psi$ is in $C^2(\RR)$, and by derivation of the equation, we show that the solution belongs to $C^\infty(\RR)$. \\
\textbf{Symmetry and positivity:} The function $\psi(-v)$ satisfies the Cauchy problem \eqref{pb de Cauchy psi}, with the same initial condition, hence we get the symmetry $\psi(-v)=\psi(v)$ on $\RR$, by uniqueness of the solution.\\
Let us now show that $\psi(v)>0$ for all $v$ in $\RR$. We check it on $\RR^+$ since $\psi$ is even. Suppose there is $v_1>0$ such that $\psi(v_1)=0$. Since $\psi(-v)=\psi(v)$ then, $\psi'(0)=0$. So, by multiplying the equation \eqref{L_0(psi)=0} by $\psi$ and integrating it over $(0,v_1)$, we obtain: $ \int_0^{v_1} \big[|\psi'|^2 + \tilde W |\psi|^2 \big]\mathrm{d}v = 0 .$
Therfore, $\psi(v)=0$ for all $v\in[0,v_0]$, which contradicts the fact that $\psi(0)=1$. Thus, $\psi(v)>0$ for all $v$ in $\RR$.\\
\textbf{Behavior at infinity:} Since $\psi$ is even, the behavior $|v|^{-\gamma}$ is excluded, because otherwise $\psi$ will be in $L^2(\RR)$, which implies thereafter that $\psi=0$. Finally, the Lemma \ref{comportement psi} gives us the existence of a unique function $R\in L^\infty([v_0,+\infty[)$ for $v_0$ large enough, with $R(v)=O( v^{ -2 })$ such that: $\psi(v)=|v|^{\gamma+1}\big(1+R(v)\big)$, up to a multiplicative constant. This constant is determined by the uniqueness of $\psi$ and it is positive since $\psi$ is positive.
\epl

\subsection{Construction of a adequate basis of solutions}
\begin{proposition}[A basis of solutions to \eqref{L_0(psi)=0}]\label{base de L_0}  
There are two positive functions $\psi_1$ and $\psi_2$, of class $C^\infty$, solutions to the equation \eqref{L_0(psi)=0} and satisfying
\begin{enumerate}
\item For all $v$ in $\RR$: $\psi_1(-v)=\psi_2(v) >0$.
\item $\{\psi_1,\psi_2\}$ forms a basis of solutions for the differential equation \eqref{L_0(psi)=0}. Moreover, the Wronskian is given by: $W\{\psi_1,\psi_2\}=\psi_1\psi_2'-\psi_1'\psi_2=1$, for all $v$ in $\RR$.
\item we have the following estimate
\begin{equation}
\psi_1(v) \lesssim \left\{\begin{array}{l}
 |v|^{-\gamma} \quad \  \mbox{for } v\geqslant v_0 ,
\\
  |v|^{1+\gamma} \quad \mbox{for } v\leqslant -v_0 .
\end{array}\right. \quad \mbox{ moreover } \quad \label{comportement de psi_1}
\psi_1(v) \sim \left\{\begin{array}{l}
 c_1|v|^{-\gamma} \quad  \ \mbox{at } +\infty ,
\\
   c_2|v|^{1+\gamma} \quad  \mbox{at } -\infty .
\end{array}\right. 
\end{equation}
Hence, the estimate and the behavior of $\psi_2$ thanks to the symmetry:
\begin{equation}\label{psi_2}
\psi_2(v) \lesssim \left\{\begin{array}{l}
 |v|^{\gamma+1} \quad  \mbox{for } v\geqslant v_0 ,
\\
  |v|^{-\gamma} \quad \ \mbox{for } v\leqslant -v_0 ,
\end{array}\right.  \quad \mbox{ and } \quad \psi_2(v) \sim \left\{\begin{array}{l}
 c_2|v|^{1+\gamma} \quad  \mbox{at } +\infty ,
\\
   c_1|v|^{-\gamma} \ \quad  \mbox{at } -\infty ,
\end{array}\right.
\end{equation}
with $v_0>0$ large enough, and where $c_1$ and $c_2$ are two given positive constants.
\end{enumerate}
\end{proposition}

\noindent \bpp \ref{base de L_0}.
Let $\psi$ be the function constructed in the previous Lemma \ref{existence de psi}. We define the two functions $\psi_1$ and $\psi_2$ by:
\begin{equation}
\psi_1(v):=\psi(v)\int_v^{+\infty}\frac{\mathrm{d}w}{\psi^2(w)} \quad \mbox{and } \ \psi_2(v)=\psi(v)\int_{-\infty}^v\frac{\mathrm{d}w}{\psi^2(w)}=\psi_1(-v) .
\end{equation}
We are going to establish the properties of $\psi_1$ and those of $\psi_2$ are obtained by symmetry.\\
The function $\psi_1$ is well defined and positive since $\psi(v)>0$ for all $v$ in $\RR$ and 
$$\frac{1}{\psi^2(v)}\underset{|v| \sim \infty}{\sim}c^{-2}|v|^{-2(\gamma+1)}\in L^1(\RR).$$ 
Moreover, 
$$ \psi_1(v) \sim \left\{\begin{array}{l}
  c_1|v|^{-\gamma} \quad \mbox{ at } +\infty ,
\\
    c_2|v|^{\gamma+1} \quad \mbox{at } -\infty .
\end{array}\right.
$$
Indeed, since $\psi(v)=c|v|^{\gamma+1}\big(1+R(v)\big)$ for $|v|\geqslant v_0$ with $R(v)=O(v^{-2})$. Then, for $|v|\geqslant v_*$ big enough $$C_1|v|^{\gamma+1} \leqslant \psi(v) \leqslant C_2|v|^{\gamma+1} .$$ 
Therefore,
$$ C_1'|v|^{\gamma+1}\int_v^{+\infty}\frac{\mathrm{d}w}{\psi^2(w)} \leqslant \psi_1(v) \leqslant C_2' |v|^{\gamma+1}\int_v^{+\infty}\frac{\mathrm{d}w}{\psi^2(w)} \ \mbox{ for all } |v|\geqslant v_*. $$
Now, for $v\geqslant v_*$
$$\int_v^{+\infty}\frac{\mathrm{d}w}{\psi^2(w)} \sim c' \int_v^{+\infty}\frac{\mathrm{d}w}{|w|^{2(\gamma+1)}} \sim c'' |v|^{-2\gamma-1} $$ 
and for $v\leqslant - v_*$
$$ \int_v^{+\infty}\frac{\mathrm{d}w}{\psi^2(w)}\leqslant \int_{-\infty}^{+\infty}\frac{\mathrm{d}w}{\psi^2(w)} = C>0 .$$
Hence the estimate \eqref{comportement de psi_1} hold.\\
The functions $\psi_1$ and $\psi_2$ are two linearly independent solutions of \eqref{L_0(psi)=0} and their Wronskian is given by:
$$ W\{\psi_1,\psi_2\}=\int_\RR \frac{\mathrm{d}w}{\psi^2(w)}=\big\|1/\psi\big\|_2^2.$$
Observe that $\|1/\psi\|_2$ is a positive constant by which we can divide in the definition of $\psi_1$ and $\psi_2$ in order to obtain $W\{\psi_1,\psi_2\}=1$.
\epp

\subsection{Solutions of $L_0(\psi)=f$ and definition of $T_0$.}
 Let's go back to the equation with a general source term \eqref{L_0(psi)=f}. We have the following Definition/Proposition:
\begin{proposition}[Definition of $T_0$]\label{def de T_0}
There exists a ``right inverse operator" of $L_0$, denoted by $T_0$,  with integral kernel $K_0$ such that, for all $f$ in $L^\infty(\RR;\langle v\rangle^{-\sigma}\mathrm{d}v)$ with $\sigma>\gamma+2$,
\begin{equation}\label{T_0}
T_0(f)(v):= \int_\RR K_0(v,w) f(w) \mathrm{d}w ,
\end{equation} 
is solution to the equation \eqref{L_0(psi)=f}, where $K_0(v,w)$ is given by:
\begin{equation}\label{K_0}
K_0(v,w):=\psi_1(v)\psi_2(w)\chi_{\{w<v\}}+\psi_1(w)\psi_2(v)\chi_{\{w>v\}} .
\end{equation}
Thus,
\begin{equation}
T_0(f)(v)= \int_{-\infty}^v \psi_1(v)\psi_2(w)f(w)\mathrm{d}w +\int_v^{\infty}\psi_1(w)\psi_2(v) f(w) \mathrm{d}w .
\end{equation}
\end{proposition}
\bp For $f$ in $L^\infty(\RR;\langle v\rangle^{-\sigma}\mathrm{d}v)$ with $\sigma>\gamma+2$, the integral \eqref{def de T_0} makes sense and $T_0$ is well defined. The function $T_0(f)$ belongs to $C^\infty(\RR)$ since $\psi_1$ and $\psi_2$ are in $C^\infty(\RR)$, in particular $T_0(f)$ is in $C^2(\RR)$. Moreover, a simple computation gives:
$$ L_0\big[T_0(f)\big]=[-\partial_v^2+\tilde W(v)]T_0(f)=f.   $$
\ep
\section{Green function for the leading part of the perturbed operator} 
We will proceed as for $L_0$, we start by studying the homogeneous equation and constructing a basis of solutions, and establishing the behavior of the latter, then we return to the equation with right order terme by defining a suitable right inverse operator.
\subsection{Approximation for large velocities}
Since $\tilde W(v) \underset{|v| \rightarrow \infty}{\longrightarrow}0$, we will approximate the homogeneous equation \begin{equation}\label{L_eta(psi)=0}
L_{\lambda,\eta}(\psi)=\big[-\pa_v^2+\tilde W(v)+ \mathrm{i} \eta v -\lambda\eta^\frac{2}{3}\big]\psi=0 
\end{equation}
by the equation 
\begin{equation}\label{tilde L_lambda,eta(psi)=0}
\tilde L_{\lambda,\eta}(\psi):=\big[-\pa_v^2+ \mathrm{i} \eta v -\lambda\eta^\frac{2}{3}\big]\psi=0  ,
\end{equation} 
which becomes after the rescaling $v=\eta^{-\frac{1}{3}}s$:
\begin{equation}\label{tilde L_lambda(phi)=0}
\tilde L_\lambda(\phi):=\big[-\pa_s^2+ \mathrm{i} s -\lambda\big]\phi=0 ,
\end{equation}
where $\phi(s):=\psi(\eta^{-\frac{1}{3}}s)$.\\

\noindent Note that the operator $\tilde L_\lambda$ is nothing but Airy's operator ``$-\partial^2_z+z$" modified. Let's give a little reminder about the Airy function.\\

\noindent \textbf{Airy's equation.} We collect the following from \cite{Lerner} and \cite{VaSo}.
\begin{lemma}[Fundamental solutions to the Airy equation]\label{Airy} 
\item Let us denote by $Ai$ the Airy function given by
\begin{equation}
Ai(z)=\frac{1}{2\pi}\int_{\mathbb{R}+\mathrm{i}\zeta}e^{\mathrm{i}(\frac{t^3}{3}+zt)}\mathrm{d}t, \ \zeta >0. 
\end{equation}
$Ai$ is independent of $\zeta > 0$ and has the following properties:
\begin{enumerate}
\item It defines an entire function of $ z \in \mathbb{C}$ and solves the Airy equation:
\begin{equation}\label{eq de Airy}
\partial^2_z Ai(z) = z Ai(z) .
\end{equation} 
\item There holds the symmetry:
$$\forall z \in \mathbb{C},\  Ai(\overline{z})=\overline{Ai(z)}.$$
\item $\ds Ai(jz)$ and $Ai(j^2z)$, where $\ds j=e^{\frac{2\pi \mathrm{i}}{3}}$, are two other fundamental solutions of \eqref{eq de Airy} with
$$ \forall z \in \mathbb{C}, \ Ai(z) + j Ai(jz) + j^2Ai(j^2) = 0 .$$
\item There holds the asympotic expansion:
\begin{equation}\label{dev.  de Airy}
Ai(z)= \left\{\begin{array}{l} \ds \frac{1}{\pi 3^{\frac{2}{3}}} \sum_{n=0}^\infty \frac{\Gamma\big(\frac{n+1}{3}\big)}{n!}\sin\bigg(\frac{2}{3}(n+1)\pi\bigg)\big(3^{\frac{1}{3}}z\big)^n ,  \ |z| \leqslant 1  ,\\
 \\
 \ds \frac{1}{2\sqrt{\pi}} z^{-\frac{1}{4}}e^{-\frac{2}{3}z^{\frac{3}{2}}} \bigg( 1+\mathcal{O}\bigg(\frac{1}{|z|^{\frac{3}{2}}}\bigg)\bigg) ,\ \forall z \in \mathbb{C}\setminus\mathbb{R}_-;  |z| \geqslant 1 ,
\end{array}\right.
\end{equation}
and simarily for derivatives.
\end{enumerate}
\end{lemma}

\noindent We will look for solutions for the equation \eqref{tilde L_lambda(phi)=0} 
under the form: $$ \phi(s) = Ai\big(a s +b\big)  .$$
We have: 
$$ \partial^2_s \phi(s) = a^2 \partial^2_s Ai(a s+b) = a^2\big(a s+ b\big)\phi(s)  .$$
Then, by identification we get
$$ a^3 = \mathrm{i} \quad \mbox{and } \ b a^2 =-\lambda .$$
Which implies that  $$ a_k = e^{\mathrm{i}\frac{\pi}{6}}e^{\mathrm{i}\frac{2\pi k}{3}}, \quad b_k = \frac{-\lambda}{a_k^2} \quad  \mbox{with } k = 0,1,2 .$$ 
Hence, each of the following three functions is a solution of \eqref{tilde L_lambda(phi)=0}:
$$ \mathfrak{a}_{\lambda}(s) = Ai\big[e^{\mathrm{i}\frac{\pi}{6}}(s+\mathrm{i}\lambda)\big], \mathfrak{b}_{\lambda}(s) = Ai\big[e^{\mathrm{i}\frac{\pi}{6}}j(s+\mathrm{i}\lambda)\big] \mbox{ and }  \mathfrak{c}_{\lambda}(s) = Ai\big[e^{\mathrm{i}\frac{\pi}{6}}j^2(s+\mathrm{i}\lambda)\big] .$$
Note that $ \big\{ \mathfrak{a}_{\lambda} , \mathfrak{b}_{\lambda}\big\} $ form a basis of solutions to the differential equation of order 2 \eqref{tilde L_lambda(phi)=0}.\\

\noindent \textbf{Back to equation \eqref{L_eta(psi)=0}.} The equation \eqref{L_eta(psi)=0} admits two possible behaviors at infinity, $\mathfrak{a}_{\lambda,\eta}$ and $\mathfrak{b}_{\lambda,\eta}$, which form a basis of solutions for the approximate equation \eqref{tilde L_lambda,eta(psi)=0}. We have the following Lemma:
\begin{lemma}[Airy behavior for large velocities for \eqref{L_eta(psi)=0}]\label{comportement de psi^eta} For $s_0>0$ large enough, there exists a unique function $R^{\lambda,\eta}_1 \in L^\infty(\{|s| \geqslant s_0\})$ (respectively, there exists a unique function $R^{\lambda,\eta}_2 \in L^\infty(\{|s| \geqslant s_0\})$) such that, the functions 
\begin{equation}\label{psi^eta rescalee}
\tilde \psi^{\lambda,\eta}_1(v)= \mathfrak{a}_{\lambda,\eta}(v)\big(1+R_1^{\lambda,\eta}(\eta^\frac{1}{3}v)\big) , \quad  |v|\geqslant s_0\eta^{-\frac{1}{3}} ,
\end{equation}
\begin{equation}\label{psi^eta rescalee 2}
\tilde \psi^{\lambda,\eta}_2(v)= \mathfrak{b}_{\lambda,\eta}(v)\big(1+R_2^{\lambda,\eta}(\eta^\frac{1}{3}v)\big) , \quad  |v|\geqslant s_0\eta^{-\frac{1}{3}}
\end{equation}
are solutions to the homogeneous equation \eqref{L_eta(psi)=0} in $\{|v|\geqslant s_0\eta^{-\frac{1}{3}}\}$, unique up to two complex multiplicative constants which can depend on $\lambda$ and $\eta$, and where $\mathfrak{a}_{\lambda,\eta}(v):=\mathfrak{a}_{\lambda}(\eta^\frac{1}{3}v)$ and $\mathfrak{b}_{\lambda,\eta}(v):=\mathfrak{b}_{\lambda}(\eta^\frac{1}{3}v)$.
Moreover, $R_i^{\lambda,\eta}$ are holomorphic in $\{|\lambda|\leqslant \lambda_0\}$ for $i=1, 2$,  with $R_1^{\lambda,\eta}(s)=O(|s|^{-\frac{3}{2}})$ and $R_2^{\lambda,\eta}(s)=O(|s|^{-\frac{3}{2}}|\mathfrak{a}_0(s)|^2)$ uniformly in $\{|\lambda|\leqslant \lambda_0\}$ and $\eta \in[0,\eta_0]$. 
\end{lemma}
\begin{remark} 
\item \begin{enumerate}
\item Thanks to the second point of Lemma \ref{Airy}, the functions $\mathfrak{a}_{\lambda,\eta}$ and $\mathfrak{b}_{\lambda,\eta}$ are linked by the following identity:
\begin{equation}\label{bar a(-v) =b(v)}
\mathfrak{\bar a}_{\bar\lambda,\eta}(-v)=\mathfrak{b}_{\lambda,\eta}(v) , \quad \forall v \in \RR .
\end{equation}
\item The Airy function $\mathfrak{a}_\lambda$ decreases exponentially in $+\infty$ and increases exponentially in $-\infty$ and inversely for $\mathfrak{b}_\lambda$:
$$ \mathfrak{a}_\lambda(s) \sim \left\{\begin{array}{l} c_1 s^{-\frac{1}{4}}e^{-\frac{\sqrt{2}}{3}s^{\frac{3}{2}}} \ \mbox{ at } +\infty  ,\\
c_2 |s|^{-\frac{1}{4}}e^{\frac{\sqrt{2}}{3}|s|^{\frac{3}{2}}}   \mbox{ at } -\infty .
\end{array}\right.
$$
\end{enumerate}
\end{remark}
\bp The proof is identical to that of Lemma \ref{comportement psi}, it relies on the fixed point theorem. \\
$\bullet$ Let's start by showing \eqref{psi^eta rescalee}. If $\tilde \psi^{\lambda,\eta}_1$ is a solution to \eqref{L_eta(psi)=0} of the form $\mathfrak{a}_{\lambda,\eta}(v)\big(1+R_1^{\lambda,\eta}(\eta^\frac{1}{3}v)\big)$ then, after rescaling the equation \eqref{L_eta(psi)=0} by $v=s\eta^{-\frac{1}{3}}$, $R^{\lambda,\eta}_1$ satisfies 
$$2\mathfrak{a}_{\lambda}'(s)(R_1^{\lambda,\eta})'(s)+\mathfrak{a}_{\lambda}(s) (R_1^{\lambda,\eta})''(s)=\tilde W_\eta(s)\mathfrak{a}_{\lambda}(s)(1+ R_1^{\lambda,\eta}(s)), $$
where $\tilde W_\eta(s) = \eta^{-\frac{2}{3}}\tilde W(\eta^{-\frac{1}{3}}s)=\frac{\gamma(\gamma+1)}{\eta^\frac{2}{3}+s^2}$. Define $\mathbb{K}_1^{\lambda,\eta}$ by
$$
\mathbb{K}_1^{\lambda,\eta}(g):=\left\{\begin{array}{l}
\ds \int_s^\infty\int_s^t \frac{\mathfrak{a}_{\lambda}^2(t)}{\mathfrak{a}_{\lambda}^2(u)}\mathrm{d}u \tilde W_\eta(t) g(t)\mathrm{d}t , \ s\geqslant s_0,
\\
\\
 \ds\int_{-\infty}^s\int_t^s \frac{\mathfrak{a}_{\lambda}^2(t)}{\mathfrak{a}_{\lambda}^2(u)}\mathrm{d}u \tilde W_\eta(t)g(t)\mathrm{d}t , \ s\leqslant -s_0  .
\end{array}\right.
$$
The previous equation on $R_1^{\lambda,\eta}$ is equivalent to $(Id-\mathbb{K}_1^{\lambda,\eta})(R_1^{\lambda,\eta})=\mathbb{K}_1^{\lambda,\eta}(1)$. \\
The computations for the fixed point for large velocities come from \cite{LebPu}, and which we will recall here for consistency. Let's start by introducing the set $U:=\{e^{\mathrm{i}\frac{\pi}{6}}(x+\mathrm{i}\lambda); x \geqslant0, |\lambda| \leqslant\lambda_0\}$. We have from \cite{QZ} or also from the equation \eqref{dev.  de Airy}, for $z\in U$; $|z| \geqslant\frac{1}{2}$:
$$ Ai(z) = \tau(z) e^{-\frac{2}{3}z^{\frac{3}{2}}} \ \mbox{ with }  \ \frac{c_0}{(1+|z|)^\frac{1}{4}} \leqslant |\tau(z)| \leqslant \frac{c_1}{(1+|z|)^\frac{1}{4}}  .$$
Let's show that for $s_0$ large enough,  
$$ 
 \big\||s|^{n+\frac{3}{2}}\mathbb{K}_1^{\lambda,\eta}(g) \big\|_{L^\infty(\{|s|\geqslant s_0\})}\leqslant \frac{1}{2} \big\||s|^n g \big\|_{L^\infty(\{|s|\geqslant s_0\})} .
$$
We are going to establish the previous inequality for $s\geqslant s_0$ and for $s \leqslant -s_0$ it is obtained in the same way, by symmetry thanks to the parity of $\tilde W_\eta$ and the identity \eqref{bar a(-v) =b(v)}. We have
$$
\big|s^{n+\frac{3}{2}}\mathbb{K}_1^{\lambda,\eta}(g)(s)\big| \leqslant \gamma(\gamma+1)s^{n+\frac{3}{2}}\int_s^\infty\int_s^t \bigg|\frac{\mathfrak{a}_{\lambda}^2(t)}{\mathfrak{a}_{\lambda}^2(u)}\bigg|\mathrm{d}u  \frac{|t^ng(t)|}{t^{n+2}} \mathrm{d}t .
$$
First, there exists a constant $C>0$ such that for all $s\geqslant 1, |\lambda|\leqslant \lambda_0$ and $t\geqslant s$, we have
$$ \int_s^t \bigg|\frac{\mathfrak{a}_{\lambda}^2(t)}{\mathfrak{a}_{\lambda}^2(u)}\bigg|\mathrm{d}u \leqslant C t^{-\frac{1}{2}} .$$
Indeed, by \eqref{dev.  de Airy}: fourth point of Lemma \ref{Airy}, and for $|\lambda|\leqslant \lambda_0$ with $\lambda_0$ small enough
$$\begin{array}{rcl}
\dps \int_s^t \bigg|\frac{\mathfrak{a}_{\lambda}^2(t)}{\mathfrak{a}_{\lambda}^2(u)}\bigg|\mathrm{d}u&\leqslant& \dps C\int_s^t e^{-\frac{4}{3}\Re([(t+\mathrm{i}\lambda)^\frac{3}{2}-(u+\mathrm{i}\lambda)^\frac{3}{2}]e^{\mathrm{i}\frac{\pi}{4}})}\mathrm{d}u\\
&=&\dps Ct\int_{\frac{s}{t}}^1e^{-\frac{4}{3}t^{\frac{3}{2}}\Re([(1+\mathrm{i}\frac{\lambda}{t})^\frac{3}{2}-(x+\mathrm{i}\frac{\lambda}{t})^\frac{3}{2}]e^{\mathrm{i}\frac{\pi}{4}})}\mathrm{d}x \\
&\leqslant & \dps Ct\int_0^\infty e^{-xt^\frac{3}{2}}\mathrm{d}x \sim Ct^{-\frac{1}{2}}\quad \mbox{if }t\geqslant s\geqslant 1 ,
\end{array}
$$
where we performed the change of variable $x=\frac{u}{t}$ in the second line and used the mean value theorem to estimate the real part of the exponent. Thus for $|\lambda|\leqslant \lambda_0$,
$$
\big|s^{n+\frac{3}{2}}\mathbb{K}_1^{\lambda,\eta}(g)(s)\big| \leqslant C\gamma(\gamma+1)s^{n+\frac{3}{2}}\int_s^\infty t^{-(n+\frac{5}{2})} \mathrm{d}t \|s^ng\|_{L^\infty(1,\infty)} .$$
Then,
$$
\|s^{n+\frac{3}{2}}\mathbb{K}_1^{\lambda,\eta}(g)\|_{L^\infty(1,\infty)}\leqslant C\frac{\gamma(\gamma+1)}{(n+\frac{3}{2})}\|s^ng\|_{L^\infty(1,\infty)}.
$$
Finally, $\mathbb{K}_1^{\lambda,\eta}$ is bounded in $L^\infty(s_0,\infty)$ with
$$
\|\mathbb{K}_1^{\lambda,\eta}\|_{\mathcal{L}(L^\infty(s_0,\infty))}\leqslant \frac{2C}{3}\gamma(\gamma+1) s_0^{-\frac{3}{2}}\leqslant \frac{1}{2}\quad \mbox{if }s_0\mbox{ is big enough.}
$$ 
The rest of the proof concerning \eqref{psi^eta rescalee} is identical to that of Lemma \ref{comportement psi}. Holomorphy of $R_1^{\lambda,\eta}$ comes from the fact that $\mathfrak{a}_\lambda$ is holomorphic on $\{|\lambda| \leqslant\lambda_0\}$, since Airy's function is an entire function in $\CC$. \\
$\bullet$ Let us study \eqref{psi^eta rescalee 2}. For a solution $\tilde \psi^{\lambda,\eta}_2$ of the form $\mathfrak{b}_{\lambda,\eta}(v)\big(1+R_2^{\lambda,\eta}(\eta^\frac{1}{3}v)\big)$, the same computations as for \eqref{psi^eta rescalee} lead to find $R^{\lambda,\eta}_2$ solution to
\begin{equation*}
\left\{\begin{array}{l} 
(Id-\mathbb{K}_2^{\lambda,\eta})(R_2^{\lambda,\eta})=\mathbb{K}_2^{\lambda,\eta}(1)  ,\\
\\
 \ds \mathbb{K}^{\lambda,\eta}_2(g) := \int_s^\infty\int_s^t \frac{\mathfrak{b}_{\lambda}^2(t)}{\mathfrak{b}_{\lambda}^2(u)}\mathrm{d}u \tilde W_\eta(t) g(t)\mathrm{d}t , \ \ s\geqslant s_0 , \\
 \\
 \ds \mathbb{K}^{\lambda,\eta}_2(g) := \int_{-\infty}^s\int_t^s \frac{\mathfrak{b}_{\lambda}^2(t)}{\mathfrak{b}_{\lambda}^2(u)}\mathrm{d}u \tilde W_\eta(t) g(t)\mathrm{d}t , \ s\leqslant -s_0 .
\end{array}\right.
\end{equation*}
As in the previous point, we establish the inequality on the norm of $\KK^{\lambda,\eta}_2$ on $[s_0,+\infty[$, and on $]-\infty,-s_0] $ it is obtained by symmetry. We have
$$
\frac{W\{\mathfrak{a}_{\lambda},\mathfrak{b}_{\lambda}\}}{\mathfrak{b}_{\lambda}^2(u)}=\bigg(\frac{\mathfrak{a}_{\lambda}(u)}{\mathfrak{b}_{\lambda}(u)}\bigg)' ,
$$
where $W\{\mathfrak{a}_{\lambda},\mathfrak{b}_{\lambda}\}$ is the Wronskian of $\mathfrak{a}_{\lambda}$ and $\mathfrak{b}_{\lambda}$ given by $W\{\mathfrak{a}_{\lambda},\mathfrak{b}_{\lambda}\}=\mathfrak{a}_{\lambda}\mathfrak{b}_{\lambda}'-\mathfrak{a}_{\lambda}'\mathfrak{b}_{\lambda}=\frac{1}{2\pi}$. Therefore,
\begin{align*}
\bigg|\frac{\mathbb{K}_2^{\lambda,\eta}(g)(s)}{s^{-(n+\frac{3}{2})}\mathfrak{a}_\lambda^2(s)}\bigg| &\leqslant C s^{n+\frac{3}{2}}\int_s^\infty\bigg(\bigg|\frac{\mathfrak{a}_{\lambda}(s)}{\mathfrak{b}_{\lambda}(s)}\bigg|+\bigg|\frac{\mathfrak{a}_{\lambda}(t)}{\mathfrak{b}_{\lambda}(t)}\bigg|\bigg)  \bigg|\frac{\mathfrak{b}_{\lambda}^2(t)}{\mathfrak{a}_{\lambda}^2(s)}\bigg| |g(t)|\frac{\mathrm{d}t}{t^2}  \\
&\leqslant C s^{n+\frac{3}{2}}\int_s^\infty\bigg(\bigg|\frac{\mathfrak{a}_{\lambda}(t)\mathfrak{b}_{\lambda}(t)}{\mathfrak{a}_{\lambda}(s)\mathfrak{b}_{\lambda}(s)}\bigg|+ \bigg|\frac{\mathfrak{a}_{\lambda}^2(t)}{\mathfrak{a}_{\lambda}^2(s)}\bigg| \bigg)|\mathfrak{a}_{\lambda}(t)\mathfrak{b}_{\lambda}(t)|   \bigg|\frac{g(t)}{t^{-n}\mathfrak{a}_{\lambda}^2(t)}\bigg| \frac{\mathrm{d}t}{t^{n+2}} ,
\end{align*}
where $C=2\pi\gamma(\gamma+1)$. Now, since we have $|Ai(z)Ai(jz)|= |\tau(z)|^2$ (note that for the following choice of determination of the complex logarithm:
$\log(re^{\mathrm{i}\theta}) = \log(r) + \mathrm{i}\theta, \ (re^{\mathrm{i}\theta})^\alpha = r^\alpha e^{\mathrm{i}\alpha\theta} $ for $\theta \in (-\pi,\pi)$,
one has $(e^{\mathrm{i}\frac{\pi}{6}}j)^{\frac{3}{2}}=-(e^{\mathrm{i}\frac{\pi}{6}})^ {\frac{3}{2}}$). Then,
for $t\geqslant s \geqslant 1$ and $|\lambda|\leqslant \lambda_0 $ with $\lambda_0$ small enough 
\begin{equation}\label{estimations a_lambda b_lambda}
|\mathfrak{a}_{\lambda}(t)\mathfrak{b}_{\lambda}(t)|=\frac{1}{4\pi}|t+\mathrm{i}\lambda|^{-\frac{1}{2}} \quad \mbox{ and }  \quad \bigg|\frac{\mathfrak{a}_{\lambda}^2(t)}{\mathfrak{a}_{\lambda}^2(s)}\bigg| \leqslant \frac{|t+\mathrm{i}\lambda|^{-\frac{1}{2}}}{|s+\mathrm{i}\lambda|^{-\frac{1}{2}}} \leqslant 2 .
\end{equation}
Then,
$$ \bigg|\frac{\mathbb{K}_2^{\lambda,\eta}(g)(s)}{s^{-(n+\frac{3}{2})}\mathfrak{a}_\lambda^2(s)}\bigg| \leqslant C \gamma(\gamma+1)s^{n+\frac{3}{2}}\int_s^\infty \frac{\mathrm{d}t}{t^{n+\frac{5}{2}}} \left\| \frac{g}{s^{-n}\mathfrak{a}_\lambda^2}\right\|_{L^\infty(s_0,\infty)} .
$$
Which implies that
$$ 
\left\|\frac{\mathbb{K}_2^{\lambda,\eta}(g)}{s^{-n}\mathfrak{a}_\lambda^2}\right\|_{L^\infty(s_0,\infty)} \leqslant C \frac{\gamma(\gamma+1)}{(n+\frac{3}{2})} s_0^{-\frac{3}{2}} \left\| \frac{g}{s^{-n}\mathfrak{a}_\lambda^2}\right\|_{L^\infty(s_0,\infty)} \leqslant \frac{1}{2}\left\| \frac{g}{s^{-n}\mathfrak{a}_\lambda^2}\right\|_{L^\infty(s_0,\infty)} 
$$
for  $s_0$ large enough. The rest of the proof is similar to the previous point, which completes the proof of the Lemma.
\ep
\subsection{Approximation for all ranges of velocities}
As for $\eta=0$, we are going to construct an intermediate function, solution to \eqref{L_eta(psi)=0}, with good properties, and which will allow us thereafter to define two other solutions forming a basis of solutions for the ODE \eqref{L_eta(psi)=0} as well as a right inverse operator for $L_{\lambda,\eta}$, which gives us solutions for the equation with right-hand side. For this, consider the following Cauchy problem:
\begin{equation}\label{pb de Cauchy psi^eta}
\left\{\begin{array}{l}
 L_{\lambda,\eta}(\psi)=0, \\
   \psi(0)=1, \\
   \psi'(0)=0.
\end{array}\right.
\end{equation}
\begin{lemma}[Existence of a global solution and properties]\label{existence psi^eta}
There exists a unique global solution $\psi^{\lambda,\eta}$ for the Cauchy problem \eqref{pb de Cauchy psi^eta}, of class $C^\infty$ with respect to $v$, holomorphic in $\{|\lambda|\leqslant\lambda_0\}$ and continuous with respect to $\eta\in[0,\eta_0]$. Furthermore,  
\begin{enumerate}
\item For all $v \in \RR$, $\bar{\psi}^{\bar\lambda,\eta}(-v)=\psi^{\lambda,\eta}(v)$ and $\psi^{\lambda,\eta}(v)\neq0$.
\item For $s_0>0$ large enough, there is a unique function $R^{\lambda,\eta} \in L^\infty[s_0,+\infty[$, holomorphic in $\{|\lambda|\leqslant \lambda_0\}$ with $R^{\lambda,\eta}(s)=O(|s|^{-\frac{3}{2}}|\mathfrak{a}_0(s)|^2)$ uniformly in $\{|\lambda|\leqslant \lambda_0\}$ and $\eta$, such that
\begin{equation}\label{psi^eta}
\psi^{\lambda,\eta}(v)= c_{\lambda}\eta^{-\frac{\gamma+1}{3}}\mathfrak{b}_{\lambda,\eta}(v)\big(1+R^{\lambda,\eta}(\eta^\frac{1}{3}v) \big) , \quad \forall v\geqslant s_0\eta^{-\frac{1}{3}} ,
\end{equation}
where $c_\lambda$ is a holomorphic function in $\{|\lambda|\leqslant\lambda_0\}$, bounded by two positive constants.
\end{enumerate}
\end{lemma}
\bp 1. \textbf{Existence:} The existence of a unique global solution $\psi^{\lambda,\eta}$ for the Cauchy problem \eqref{pb de Cauchy psi^eta} follows from the Cauchy-Lipschitz theorem. \\
\textbf{Symmetry:} For all $v$ in $\RR$, ${\bar \psi}^{\bar \lambda,\eta}(-v)$ satisfies the equation \eqref{pb de Cauchy psi^eta}. Hence, ${\bar \psi}^{\bar \lambda,\eta}(-v)=\psi^{\lambda,\eta}(v)$ by uniqueness of the solution. Let us now show that the solution $\psi^{\lambda,\eta}$ does not vanish on $\RR$, for all $\lambda \in \CC; |\lambda| \leqslant\lambda_0$ and $\eta>0$. We verify it only on $\RR^+$ thanks to symmetry. The idea of the proof was inspired from \cite{LebPu} to prove that the rescaled function does not vanish on $\RR^+$. Suppose there is $v_1>0$ such that $\psi^{\lambda,\eta}(v_1)=0$. Set $v_1=s_1\eta^{-\frac{1}{3}}$. By integrating the equation of $\psi^{\lambda,\eta}$ rescaled by $v=s\eta^{-\frac{1}{3}}$, and multiplied by $\Psi_{\lambda,\eta}(s):=\psi^{\lambda ,\eta}(s\eta^{-\frac{1}{3}})$, between $0$ and $s_1$, we obtain:
$$ \int_0^{s_1} \bigg[|(\Psi_{\lambda,\eta})'|^2 + \frac{\gamma(\gamma+1)}{\eta^\frac{2}{3}+s^2}|\Psi_{\lambda,\eta}|^2 \mathrm{d}s+ \mathrm{i} s|\Psi_{\lambda,\eta}|^2\bigg]\mathrm{d}s=\lambda\int_0^{s_1} |\Psi_{\lambda,\eta}|^2\mathrm{d}s . $$
Which implies that 
$$ \Re \lambda \geqslant \frac{\gamma(\gamma+1)}{\eta^\frac{2}{3}+s^2} \ \mbox{ and } \  \Imm \lambda \geqslant s \geqslant0 .$$
Then,  
$$ \Re \lambda + \Imm \lambda \geqslant \underset{s\in[0,s_1]}{\min} \bigg(s+ \frac{\gamma(\gamma+1)}{\eta^\frac{2}{3}+s^2}\bigg) \geqslant c_0 > 0  \ \mbox{ for } \eta \in[0,\eta_0] .$$
Which contradicts the fact that $|\lambda|\leqslant \lambda_0$, since $c_0$ does not depend on $\lambda_0$ and we can choose $\lambda_0$ as small as we want.  \\
2. \textbf{Airy behavior:} Since there are only two possible behaviors at infinity \cite{LebPu}: $\mathfrak{a}_{\lambda,\eta}$ and $\mathfrak{b}_{\lambda,\eta}$, using exactly the same argument  that we used to show that $\psi^{\lambda,\eta}\neq0$, we can show that the profile $\mathfrak{a}_{\lambda,\eta}$, Airy which decreases exponentially at $+\infty$, so in $L^2(1,\infty)$, is excluded. Indeed, by rescaling the equation by $v=\eta^{-\frac{1}{3}}s$ and by multiplying it by $\Psi_{\lambda,\eta}(s):=\psi^{\lambda,\eta}(\eta^{-\frac{1}{3}}s)$ and integrating it over $(0,\infty)$, thanks to the condition $(\psi^{\lambda,\eta})'(0)=0$ (recalling that $\psi^{\lambda,\eta}$ satisfies the Cauchy problem \eqref{pb de Cauchy psi^eta}), we obtain:
$$\int_0^\infty \bigg[|(\Psi_{\lambda,\eta})'|^2 + \frac{\gamma(\gamma+1)}{\eta^\frac{2}{3}+s^2}|\Psi_{\lambda,\eta}|^2 \mathrm{d}s+\mathrm{i} s|\Psi_{\lambda,\eta}|^2\bigg]\mathrm{d}s=\lambda\int_0^\infty |\Psi_{\lambda,\eta}|^2\mathrm{d}s . $$
Then,  
$$ \Re \lambda + \Imm \lambda \geqslant \underset{s\geqslant 0}{\min} \bigg(s+ \frac{\gamma(\gamma+1)}{\eta^\frac{2}{3}+s^2}\bigg) \geqslant c_0 > 0  \ \mbox{ for } \eta \in[0,\eta_0] .$$
Hence the contradiction with $|\lambda|\leqslant \lambda_0$, with $\lambda_0$ small enough. Thus, the function $\psi^{\lambda,\eta}$ takes the profile of Airy which explodes at $+\infty$, and therefore $\psi^{\lambda,\eta}(v) \underset{ +\infty}{\sim} \mathfrak{b}_{\lambda,\eta}(v)$. To obtain the equation \eqref{psi^eta}, by the previous Lemma, there exists a unique function $R^{\lambda,\eta} \in L^\infty(\{|s| \geqslant s_0\})$ such that $\psi^{\lambda,\eta}$ takes the form \eqref {psi^eta rescalee 2}, up to a multiplicative constant. This constant is determined by the uniqueness of $\psi^{\lambda,\eta}$. To calculate it, we use the continuity of $\psi^{\lambda,\eta}$ with respect to $v$ at the point $s_0 \eta^{-\frac{1}{3}}$ and the fact that $\psi^{\lambda,\eta}$ is close to $\psi$ on $[-s_0 \eta^{-\frac{1}{3}},s_0 \eta^{-\frac{1} {3}}]$. For this last point, we will prove it in the following lemma.
\ep
\begin{lemma}[Some estimates on $\psi^{\lambda,\eta}$]\label{estimation psi^eta}
\item \begin{enumerate}
\item For all $\sigma \in (0,\frac{1}{3})$, for all $a>0$ we have the following uniform convergence: 
\begin{equation}\label{psi^eta-->psi}
\left\|\frac{(\psi^{\lambda,\eta}-\psi)}{\psi}\right\|_{L^\infty \big(-a\eta^{\sigma-\frac{1}{3}},a\eta^{\sigma-\frac{1}{3}}\big)} \underset{\eta \rightarrow  0}{\longrightarrow} 0.
\end{equation}
\item There exist two positive constants $C_1$ and $C_2$ such that,  the following estimates holds:\\
$\bullet$ For all $|v|\leqslant s_0\eta^{-\frac{1}{3}}$, we have 
\begin{equation}\label{psi^eta bornee par psi}
C_1 \psi(v) \leqslant |\psi^{\lambda,\eta}(v)| \leqslant C_2 \psi(v) .
\end{equation}
$\bullet$ For all $v$ in $\RR$, we have
\begin{equation}\label{psi^eta>psi}
\psi(v) \leqslant C_2 |\psi^{\lambda,\eta}(v)| .
\end{equation}
\item We have the limit 
\begin{equation}\label{1/(psi^eta)2-->1/(psi)2}
\int_\RR \frac{\mathrm{d}w}{[\psi^{\lambda,\eta}(w)]^2} \ \underset{\eta \rightarrow  0}{\longrightarrow} \ \int_\RR \frac{\mathrm{d}w}{[\psi(w)]^2}  .
\end{equation}
\end{enumerate}
\end{lemma}

\noindent \bp 1. Let $\sigma\in(0,\frac{1}{3})$ and $a>0$. Thanks to the symmetry, we show the limit \eqref{psi^eta-->psi} on $[0,a\eta^{\sigma-\frac{1}{3}}]$. Let's set $\phi^{\lambda,\eta}:=\psi^{\lambda,\eta}-\psi$. Then, $\phi^{\lambda,\eta}$ satisfies the equation:
$$ L_0[\phi^{\lambda,\eta}](v)=[-\pa_v^2+\tilde{W}(v)]\phi^{\lambda,\eta}(v)=(\lambda\eta^{\frac{2}{3}}- \mathrm{i} \eta v)[\phi^{\lambda,\eta}(v)+\psi(v)]$$
First, recall that $\psi^{\lambda,\eta}(0)=\psi(0)=1$ and $(\psi^{\lambda,\eta})'(0)=\psi'(0)=0$. Thus, $\phi^{\lambda,\eta}$ satisfies the initial conditions $\phi^{\lambda,\eta}(0)=(\psi^{\lambda,\eta})'(0)= 0$. The solution of the equation $L_0(\phi)=S$ with $\phi(0)=\phi'(0)=0$ is given by:
$$ \phi(v)=\psi_1(v)\int_0^v \psi_2(w)S(w)\mathrm{d}w-\psi_2(v)\int_0^v \psi_1(w)S(w)\mathrm{d}w $$
 taking into account the fact that $\psi_1(-v)=\psi_2(v)$. Therefore, $\phi^{\lambda,\eta}/\psi$ satisfies:
\begin{align*}
\bigg|\frac{\phi^{\lambda,\eta}(v)}{\psi(v)}\bigg| &\leqslant \frac{1}{\psi(v)}\int_0^v \big[\psi_1(v)\psi_2(w)+\psi_1(w)\psi_2(v)\big]\big|(\lambda\eta^{\frac{2}{3}}-\mathrm{i} \eta w)[\phi^{\lambda,\eta}(w)+\psi(w)]\big|\mathrm{d}w \\
&\leqslant \frac{2}{\psi(v)} \bigg(\int_0^v(|\lambda|\eta^{\frac{2}{3}}+\eta |w|)\big[\psi_1(v)\psi_2(w)+\psi_1(w)\psi_2(v)\big]\psi(w)\mathrm{d}w\bigg) \bigg(1+\left\|\frac{\phi^{\lambda,\eta}}{\psi}\right\|_\infty\bigg) .
\end{align*}
Now we split the two cases, $v\in[0,v_0]$ and $v\in[v_0,a\eta^{\sigma-\frac{1}{3}}]$ where $v_0>0$ is large enough such that $\psi_1$ and $\psi_2$ satisfy \eqref{comportement de psi_1} and \eqref{psi_2}.

\noindent \textbf{Step 1: $v\in[0,v_0]$}. We have in this case, $\psi$, $\psi_1$ and $\psi_2$ are bounded above and below by positive constants. Therefore,
\begin{align*}
\bigg|\frac{\phi^{\lambda,\eta}(v)}{\psi(v)}\bigg| &\lesssim \bigg(\int_0^{v_0}(|\lambda|\eta^{\frac{2}{3}}+\eta |w|)\big[\psi_2(w)+\psi_1(w)\big]\psi(w)\mathrm{d}w\bigg) \bigg(1+\left\|\frac{\phi^{\lambda,\eta}}{\psi}\right\|_\infty\bigg) \\ 
&\lesssim \big(|\lambda|\eta^{\frac{2}{3}}+\eta\big)\bigg(1+\left\|\frac{\phi^{\lambda,\eta}}{\psi}\right\|_\infty\bigg) .
\end{align*}
Therefore, for $\eta$ small enough 
\begin{equation}\label{estim.  phi^eta/psi v petit}
\left\|\frac{\phi^{\lambda,\eta}}{\psi}\right\|_\infty \lesssim |\lambda|\eta^{\frac{2}{3}}+\eta .
\end{equation}
\textbf{Step 2: $v\in[v_0,a\eta^{\sigma-\frac{1}{3}}]$}. We have in this case, $\psi(v)\sim v^{\gamma+1} \sim \psi_2(v)$ and $\psi_1(v) \sim v^{-\gamma}$, up to a multiplicative constants, and since these three functions are bounded on the compact $[0,v_0]$ then,
\begin{align*}
\bigg|\frac{\phi^{\lambda,\eta}(v)}{\psi(v)}\bigg| &\lesssim \bigg(\bigg|\frac{\psi_1(v)}{\psi(v)}\bigg|\int_0^v(|\lambda|\eta^{\frac{2}{3}}+\eta |w|)\psi_2(w)\psi(w)\mathrm{d}w \\
&+\bigg|\frac{\psi_2(v)}{\psi(v)}\bigg|\int_0^v(|\lambda|\eta^{\frac{2}{3}}+\eta |w|)\psi_1(w)\psi(w)\mathrm{d}w\bigg) \bigg(1+\left\|\frac{\phi^{\lambda,\eta}}{\psi}\right\|_\infty\bigg)\\
&\lesssim \big[(|\lambda|\eta^{\frac{2}{3}}+\eta ) + (|\lambda|\eta^{\frac{2}{3}}v^2+\eta v^3)\big]\bigg(1+\left\|\frac{\phi^{\lambda,\eta}}{\psi}\right\|_\infty\bigg) .
\end{align*}
Thus,
\begin{equation}\label{estim.  phi^eta/psi v grand}
\left\|\frac{\phi^{\lambda,\eta}}{\psi}\right\|_\infty \lesssim a^2|\lambda|\eta^{2\sigma}+a^3\eta^{3\sigma} , \quad \forall v\in[v_0,a\eta^{\sigma-\frac{1}{3}}] .
\end{equation}
From where, with \eqref{estim.  phi^eta/psi v petit}, we get the limit \eqref{psi^eta-->psi}. \\
2. From \eqref{estim.  phi^eta/psi v petit} and \eqref{estim.  phi^eta/psi v grand}, and for $\sigma =0$ and $a$ small enough such that, for all $v$ in $[0,a\eta^{-\frac{1}{3}}]$
$$ \bigg|\frac{\phi^{\lambda,\eta}(v)}{\psi(v)}\bigg|  \leqslant \frac{1}{2}\bigg(1+\left\|\frac{\phi^{\lambda,\eta}}{\psi}\right\|_\infty\bigg)$$
we get 
$$\left\|\frac{\phi^{\lambda,\eta}}{\psi}\right\|_\infty \leqslant 1 , \quad \forall v \in [0,a\eta^{-\frac{1}{3}}] ,$$
Then, 
\begin{equation}\label{c_1 psi < psi^eta < c_2 psi}
\psi(v) \leqslant |\psi^{\lambda,\eta}(v)|\leqslant 2 \psi(v) , \quad \forall v \in [0,a\eta^{-\frac{1}{3}}] .
\end{equation}
It remains to establish the inequality \eqref{psi^eta bornee par psi} on $[a\eta^{-\frac{1}{3}},s_0\eta^{-\frac{1}{3}}]$, which is obtained by rescaling the equation of $\psi^{\lambda,\eta}$ by $v=\eta^{-\frac{1}{3}}s$.
Indeed, if we denote $\Psi_{\lambda,\eta}(s) :=\psi^{\lambda,\eta}(\eta^{-\frac{1}{3}}s)$, then $\Psi_{\lambda,\eta}$ satisfies
$$ \big[-\pa_s^2 + \frac{\gamma(\gamma+1)}{\eta^{\frac{2}{3}}+s^2} + \mathrm{i} s - \lambda \big]\Psi_{\lambda,\eta}(s) = 0 . $$
We have by a fixed point argument as in Lemma \ref{existence psi^eta}, $\Psi_{\lambda,\eta}(s)=\Psi_{\lambda,0}(s)(1+r_{\lambda, \eta}(s))$ up to a multiplicative constant, which may depend on $\lambda$ and $\eta$, and where $\Psi_{\lambda,0}$ is the solution of the last equation with $\eta=0$: 
$$ \big[-\pa_s^2 + \frac{\gamma(\gamma+1)}{s^2} + \mathrm{i} s - \lambda \big]\Psi_{\lambda,0}(s) = 0 ,$$
with $r_{\lambda,\eta}$ tends to $0 $ when $\eta$ tends to $0$.
Now, since the function $\Psi_{\lambda,0}$ is continuous on the compact $[a,s_0]$ and holomorphic in $\{|\lambda| \leqslant \lambda_0 \}$ then, since $\psi_{\lambda,\eta}$ does not vanish, $|\Psi_{\lambda,0}|$ is bounded from below and above by two positive constants, uniformly with respect to $\lambda$ and $\eta$. Then, since $r_{\lambda,\eta}$ is also bounded on $ [a,s_0]$, uniformly with respect to $\lambda$ and $\eta$ then, by continuity of $\psi^{\lambda,\eta}$ at $a\eta^{-\frac{1}{3}}$ and the inequality \eqref{c_1 psi < psi^eta < c_2 psi}, we get
$\psi^{\lambda,\eta}(v)=c_{\lambda}\eta^{-\frac{\gamma+1}{3}}\Psi_{\lambda,0}(s)(1+r_{\lambda,\eta}(s))$ with $c_\lambda$ holomorphic in $\{|\lambda| \leqslant \lambda_0 \}$. Thus, $\psi^{-1}\psi^{\lambda,\eta}\sim\eta^{\frac{\gamma+1}{3}}\psi^{\lambda,\eta}$ is bounded on $[a\eta^{-\frac{1}{3}},s_0\eta^{-\frac{1}{3}}]$ from below and above.
Hence the inequality \eqref{psi^eta bornee par psi} holds on $ [-s_0\eta^{-\frac{1}{3}},s_0\eta^{-\frac{1}{3}}]$.\\
For the inequality \eqref{psi^eta>psi}, it comes from \eqref{psi^eta bornee par psi} for $|v|\leqslant s_0\eta^{-\frac{1}{3}}$ and the fact that for $v\geqslant s_0\eta^{-\frac{1}{3}}$,  $1\lesssim |\mathfrak{b}_{\lambda,\eta}(v)|$ and $\psi(v)\lesssim \eta^{-\frac{\gamma+1}{3}}$. The case $v\leqslant -s_0\eta^{-\frac{1}{3}}$ is obtained by symmetry.
\\
3.  This limit is a direct consequence of the second point, inequality \eqref{psi^eta>psi}, and Lebesgue's theorem, since $\psi^{\lambda,\eta}$ is continuous with respect to $\eta$, which gives us simple convergence.
\ep
\begin{proposition}[Basis of solutions and definition of $T_{\lambda,\eta}$]\label{base de L_eta}
\item $\bullet$ There are two functions $\psi^{\lambda,\eta}_1$ and $\psi^{\lambda,\eta}_2$, solutions to the equation \eqref{L_eta(psi)=0}, belonging to $C^\infty( \RR,\CC)$, continuous with respect to $\eta\in[0,\eta_0]$ and holomorphic in $\{|\lambda|\leqslant \lambda_0\}$. Furthermore,
\begin{enumerate}
\item $\{\psi^{\lambda,\eta}_1,\psi^{\lambda,\eta}_2\}$ forms a basis of solutions for the equation \eqref{L_eta(psi)=0} with $W\{\psi^ {\lambda,\eta}_1,\psi^{\lambda,\eta}_2\}=1$.
\item For all $v$ in $\RR$, $\bar \psi^{\bar\lambda,\eta}_1(-v)=\psi^{\lambda,\eta}_2(v)$.
\end{enumerate}
$\bullet$ The operator $T_{\lambda,\eta}$ defined by the integral kernel $K_{\lambda,\eta}$, is a right inverse of $L_{\lambda,\eta}$:
\begin{equation}\label{T_eta(f)}
T_{\lambda,\eta}(f)(v)= \int_\RR K_{\lambda,\eta}(v,w) f(w) \mathrm{d}w  ,
\end{equation}
with
\begin{equation}\label{K_eta}
K_{\lambda,\eta}(v,w)=\psi_1^{\lambda,\eta}(v)\psi_2^{\lambda,\eta}(w)\chi_{\{w<v\}}+\psi_1^{\lambda,\eta}(w)\psi_2^{\lambda,\eta}(v)\chi_{\{w>v\}} ,
\end{equation}
continuous with respect to $(v, w) \in \RR\times\RR$ and $\eta \in [0, \eta_0]$ and holomorphic in $\{|\lambda| \leqslant \lambda_0\}$. \\
Thus, $T_{\lambda,\eta}(f)$ is solution to the equation $L_{\lambda,\eta}(\psi)=f$.
\end{proposition}
\bp $\bullet$ Let $\psi^{\lambda,\eta}$ be the function of lemma \ref{existence psi^eta}. We define the two functions $\psi^{\lambda,\eta}_1$ and $\psi^{\lambda,\eta}_2$, for $v \in \RR$, by:
\begin{equation}
\psi^{\lambda,\eta}_1(v):=\frac{\psi^{\lambda,\eta}(v)}{\|1/\psi^{\lambda,\eta}\|_2}\int_v^{+\infty}\frac{\mathrm{d}w}{[\psi^{\lambda,\eta}(w)]^2} \quad \mbox{and } \ \psi^{\lambda,\eta}_2(v)=\frac{\psi^{\lambda,\eta}(v)}{\|1/\psi^{\lambda,\eta}\|_2}\int_{-\infty}^v\frac{\mathrm{d}w}{[\psi^{\lambda,\eta}(w)]^2} .
\end{equation}
The functions $\psi^{\lambda,\eta}_1$ and $\psi^{\lambda,\eta}_2$ are both well defined thanks to the inequality \eqref{psi^eta>psi} which guarantees that the integral is indeed defined in both cases, and they are solutions to \eqref{L_eta(psi)=0} having the same regularity as $\psi^ {\lambda, \eta}$, i.e. $\psi^{\lambda,\eta}_i \in C^\infty(\RR)$ for $i=1,2$. The continuity/holomorphy of $\psi^{\lambda,\eta}_i$ for $i=1,2$ is obtained by Lebesgue's theorem thanks to the continuity/holomorphy of $\psi^{\lambda,\eta}$, the limit \eqref{1/(psi^eta)2-->1/(psi)2} and the domination $\frac{1}{|\psi^{\lambda,\eta}(v)|^2}\lesssim \frac{1}{\psi^2(v)}$ by \eqref{psi^eta>psi}. Moreover,\\
1.  We have: $W\{\psi^{\lambda,\eta}_1,\psi^{\lambda,\eta}_2\}=\psi^{\lambda,\eta}_1(\psi^{\lambda,\eta}_2)'-(\psi^{\lambda,\eta}_1)'\psi^{\lambda,\eta}_2=\|1/\psi^{\lambda,\eta}\|_2^{-2}\int_\RR \frac{\mathrm{d}w}{[\psi^{\lambda,\eta}(w)]^2}=1$.   Thus, $\psi^{\lambda,\eta}_1$ and $\psi^{\lambda,\eta}_2$ are linearly independent. \\
2. The second point comes from the symmetry of $\psi^{\lambda,\eta}$ (first point of the lemma \ref{existence psi^eta}).\\
$\bullet$ For $f$ in $L^\infty\big(\RR;\frac{\langle v\rangle^{-\sigma} \mathrm{d}v}{|\psi^{\lambda,\eta}_1|+|\psi^{\lambda,\eta}_2|}\big)$ with $\sigma>2$, the integral \eqref{T_eta(f)} is well defined, $T_{\lambda,\eta}(f)$ belongs to $C^2(\RR,\CC)$ and we have: $L_{\lambda,\eta}[T_{\lambda,\eta}(f)]=f$.
\ep

\ni \begin{remark}\label{remarques sur T_eta}
\item \begin{enumerate}
\item Note that $T_{\lambda,0}=:T_0$ does not depend on $\lambda$, since the functions $\psi_i^{\lambda,\eta}$ for $i=1,2$ are continuous with respect to $\eta$ and their limits, when $\eta\rightarrow 0$, $\psi_i$ do not depend on $\lambda$.
\item Thanks to the identity $\bar \psi^{\bar\lambda,\eta}_1(-v)=\psi^{\lambda,\eta}_2(v)$, $K_{\lambda,\eta}$ satisfies
$$ \bar K_{\bar\lambda,\eta}(-v,-w)=K_{\lambda,\eta}(v,w), \quad \forall v, w \in \RR .$$
Thus, for $\bar f_{\bar\lambda,\eta}(-v)= f_{\lambda,\eta}(v)$, we get the following identity on $T_{\lambda,\eta}$
$$ \bar T_{\bar \lambda,\eta} [f](-v) = T_{\lambda,\eta}[f](v) , \quad \forall v \in \RR  .  $$
\end{enumerate}
\end{remark}
\section{Properties of the Green functions corresponding to the leading part of the operator}
\subsection{Some complementary estimates of the the elements of the basis}
\begin{lemma}[Some estimates on $\psi^{\lambda,\eta}_1$ and $\psi^{\lambda,\eta}_2$]\label{estimation des psi^eta_i}
\item \begin{enumerate}
\item We have the following estimates:
\begin{equation}\label{psi_1^eta}
 |\psi_1^{\lambda\eta}(v)| \lesssim \left\{\begin{array}{l}
 \eta^{\frac{\gamma}{3}} |\mathfrak{a}_{\lambda,\eta}(v)| ,\qquad \ v\geqslant s_0\eta^{-\frac{1}{3}} ,
\\
  \psi_1(v) , \qquad  \qquad \ |v|\leqslant s_0\eta^{-\frac{1}{3}} ,\\
 \eta^{-\frac{\gamma+1}{3}} |\mathfrak{a}_{\lambda,\eta}(v)| , \quad  v\leqslant -s_0\eta^{-\frac{1}{3}}  .
\end{array}\right. 
\end{equation}
Similarly for $\psi^{\lambda,\eta}_2$,
\begin{equation}\label{psi_2^eta}
 |\psi_2^{\lambda\eta}(v)| \lesssim \left\{\begin{array}{l}
 \eta^{-\frac{\gamma+1}{3}} |\mathfrak{b}_{\lambda,\eta}(v)| , \quad v\geqslant s_0\eta^{-\frac{1}{3}}  ,
\\
  \psi_2(v) , \qquad \qquad \ |v|\leqslant s_0\eta^{-\frac{1}{3}} , \\
 \eta^{\frac{\gamma}{3}} |\mathfrak{b}_{\lambda,\eta}(v)| , \qquad \ v\leqslant -s_0\eta^{-\frac{1}{3}}  .
\end{array}\right. 
\end{equation}
\item There is a constant $C>0$ such that, for all $\sigma \in (0,\frac{1}{3})$, $a>0$ and for all $v \in [-a\eta^{-\frac{1}{3}},a\eta^{-\frac{1}{3}}]$ one has
\begin{equation}\label{estimation de psi^eta_i - psi_i}
\bigg| \frac{\psi^{\lambda,\eta}_i(v)-\psi_i(v)}{\psi_i(v)}\bigg| \leqslant C (a^2|\lambda| \eta^{2\sigma}+a^3\eta^{3\sigma}) , \qquad \forall i=1,2   .
\end{equation}
\end{enumerate}
\end{lemma}
\bp 1. We are going to show the inequality \eqref{psi_1^eta}, and that on $\psi^{\lambda,\eta}_2$ is obtained by symmetry. \\
\noindent $\bullet$ \textbf{Case 1: $v\geqslant s_0 \eta^{-\frac{1}{3}}$.} We have $C_1\leqslant \big|1+R^{\lambda,\eta}(\eta^\frac{1}{3}v)\big|\leqslant C_2 $ since $R^{\lambda,\eta}(s)=O(|s|^{-\frac{3}{2}}|\mathfrak{a}_0(|s|)|^2)$ for $|s|\geqslant s_0$ with $s_0$ big enough. Then,
$$
\bigg| \frac{\psi^{\lambda,\eta}_1(v)}{\eta^{\frac{\gamma}{3}} \mathfrak{a}_{\lambda,\eta}(v)} \bigg| \lesssim \eta^{-\frac{2\gamma+1}{3}}\bigg|\frac{\mathfrak{b}_{\lambda,\eta}(v)}{\mathfrak{a}_{\lambda,\eta}(v)}\bigg| \int_v^\infty \frac{\mathrm{d}w}{|\mathfrak{b}_{\lambda,\eta}(w)|^2}\eta^{\frac{2(\gamma+1)}{3}} .
$$
By performing the changes of variables $v=s\eta^{-\frac{1}{3}}$ and $w=t\eta^{-\frac{1}{3}}$,  and for $|\lambda|\leqslant \lambda_0$ with $\lambda_0$ small enough, we obtain 
\begin{align*}
\bigg| \frac{\psi^{\lambda,\eta}_1(v)}{\eta^{\frac{\gamma}{3}} \mathfrak{a}_{\lambda,\eta}(v)} \bigg| &\lesssim \bigg|\frac{\mathfrak{b}_\lambda(s)}{\mathfrak{a}_\lambda(s)}\bigg| \int_s^\infty \frac{\mathrm{d}t}{|\mathfrak{b}_\lambda(t)|^2} \\
&= 4\pi e^{\frac{4}{3}\Re[e^{\mathrm{i}\frac{\pi}{4}}(s+\mathrm{i}\lambda)^{\frac{3}{2}}]}\int_s^\infty |t+\mathrm{i}\lambda|^{\frac{1}{2}}e^{-\frac{4}{3}\Re[e^{\mathrm{i}\frac{\pi}{4}}(t+\mathrm{i}\lambda)^{\frac{3}{2}}]}\mathrm{d}t \\
&= e^{\frac{4}{3}\Re[e^{\mathrm{i}\frac{\pi}{4}}(s+\mathrm{i}\lambda)^{\frac{3}{2}}]} \int_s^\infty \Re[e^{\mathrm{i}\frac{\pi}{4}}(t+\mathrm{i}\lambda)^{\frac{1}{2}}]\frac{|t+\mathrm{i}\lambda|^{\frac{1}{2}}}{\Re[e^{\mathrm{i}\frac{\pi}{4}}(t+\mathrm{i}\lambda)^{\frac{1}{2}}]} e^{-\frac{4}{3}\Re[e^{\mathrm{i}\frac{\pi}{4}}(t+\mathrm{i}\lambda)^{\frac{3}{2}}]}\mathrm{d}t  \\
&\lesssim e^{\frac{4}{3}\Re[e^{\mathrm{i}\frac{\pi}{4}}(s+\mathrm{i}\lambda)^{\frac{3}{2}}]} \int_s^\infty \Re[e^{\mathrm{i}\frac{\pi}{4}}(t+\mathrm{i}\lambda)^{\frac{1}{2}}] e^{-\frac{4}{3}\Re[e^{\mathrm{i}\frac{\pi}{4}}(t+\mathrm{i}\lambda)^{\frac{3}{2}}]}\mathrm{d}t  \\
&\lesssim 1
\end{align*}
since for $t\geqslant s \geqslant s_0$ and $|\lambda| \leqslant \lambda_0$:
$$ \frac{|t+\mathrm{i}\lambda|^{\frac{1}{2}}}{\Re[e^{\mathrm{i}\frac{\pi}{4}}(t+\mathrm{i}\lambda)^{\frac{1}{2}}]} \lesssim 1 \ \mbox{ and } \  \int_s^\infty \Re[e^{\mathrm{i}\frac{\pi}{4}}(t+\mathrm{i}\lambda)^{\frac{1}{2}}] e^{-\frac{4}{3}\Re[e^{\mathrm{i}\frac{\pi}{4}}(t+\mathrm{i}\lambda)^{\frac{3}{2}}]}\mathrm{d}t = \frac{3}{4} e^{-\frac{4}{3}\Re[e^{\mathrm{i}\frac{\pi}{4}}(s+\mathrm{i}\lambda)^{\frac{3}{2}}]} .$$
$\bullet$ \textbf{Case 2: $v\leqslant -s_0 \eta^{-\frac{1}{3}}$.} As in the previous case, $|1+R^{\lambda,\eta}(\eta^{\frac{1}{3}}|v|)|$ is bounded below and above and since $\psi^{\lambda,\eta}(v)=\bar c_{\bar \lambda} \eta^{-\frac{\gamma+1}{3}}\mathfrak{a}_{\lambda,\eta}(v)\big(1+\bar R^{\bar\lambda,\eta}(-\eta^{\frac{1}{3}}v)\big)$ then, by \eqref{psi^eta-->psi}
$$
\bigg| \frac{\psi^{\lambda,\eta}_1(v)}{\eta^{-\frac{\gamma+1}{3}} \mathfrak{a}_{\lambda,\eta}(v)} \bigg| \lesssim \int_v^\infty \frac{\mathrm{d}w}{|\psi^{\lambda,\eta}(w)|^2} \lesssim \int_v^\infty \frac{\mathrm{d}w}{|\psi(w)|^2} \lesssim \int_\RR \frac{\mathrm{d}w}{|\psi(w)|^2} \lesssim 1 .
$$
$\bullet$ \textbf{Case 3: $|v|\leqslant s_0 \eta^{-\frac{1}{3}}$.} We have in this case by \eqref{psi^eta bornee par psi}, $ C_1 \psi(v) \leqslant |\psi^{\lambda,\eta}(v)|\leqslant C_2 \psi(v)$ and thanks to \eqref{psi^eta>psi} we get
$$ \bigg|\frac{\psi^{\lambda,\eta}_1(v)}{\psi_1(v)}\bigg| =  \bigg|\frac{\psi^{\lambda,\eta}(v)}{\psi(v)} \bigg(\int_v^\infty \frac{\mathrm{d}w}{|\psi^{\lambda,\eta}(w)|^2}\bigg) \bigg(\int_v^\infty \frac{\mathrm{d}w}{|\psi(w)|^2}\bigg)^{-1}\bigg| \lesssim 1   .  $$
2. The proof of this point is the same as for \eqref{psi^eta-->psi}.

\ep
We will end this subsection with a lemma on the estimation of the kernel $K_{\lambda,\eta}$.

\begin{lemma}
Let $\sigma \in (0,\frac{1}{3})$ and let $a>0$. The following estimate
\begin{equation}\label{estimation de K_eta - K_0}
|K_{\lambda,\eta}(v,w)-K_0(v,w)| \lesssim \big(a^2|\lambda|\eta^{2\sigma}+a^3\eta^{3\sigma}\big)K_0(v,w)
\end{equation}
holds for $|v|, |w| \leqslant a\eta^{\sigma-\frac{1}{3}}$. Therefore, for $\eta$ small enough or for $\sigma=0$ and $a$ small enough, we have the estimate
\begin{equation}\label{K_eta < K_0}
|K_{\lambda,\eta}(v,w)| \lesssim K_0(v,w) ,  \qquad \forall |v|, |w| \leqslant a\eta^{\sigma-\frac{1}{3}}.
\end{equation}
\end{lemma}
\bp Let $\sigma \in (0,\frac{1}{3})$ and $a>0$.  Denote $\phi^{\lambda,\eta}_i := \psi^{\lambda,\eta}_i-\psi_i$ for $i=1,2$. We have
\begin{align*}
K_{\lambda,\eta}(v,w)-K_0(v,w)&=\phi^{\lambda,\eta}_1(v)\psi_2^{\lambda,\eta}(w)\chi_{\{w<v\}}+\phi^{\lambda,\eta}_1(w)\psi_2^{\lambda,\eta}(v)\chi_{\{w>v\}} \\
&+ \psi_1(v) \phi^{\lambda,\eta}_2(w)\chi_{\{w<v\}}  + \psi_1(w) \phi^{\lambda,\eta}_2(v)\chi_{\{w>v\}} .
\end{align*}
By \eqref{estimation de psi^eta_i - psi_i},  $|\phi^{\lambda,\eta}_i(v)| \lesssim \big(a^2|\lambda|\eta^{2\sigma}+a^3\eta^{3\sigma}\big)\psi_i(v)$ and by \eqref{psi_1^eta} and \eqref{psi_2^eta}, $|\psi^{\lambda,\eta}_i(v)| \lesssim \psi_i(v)$ for all $|v|\leqslant a\eta^{\sigma-\frac{1}{3}}$ and $i=1,2$.  Hence inequality \eqref{estimation de K_eta - K_0} holds and for $\eta$ small enough or $a$ small enough with $\sigma\in[0,\frac{1}{3}[$ we get estimate \eqref{K_eta < K_0}.
\ep
\subsection{Continuity of the functional, introduction of the weighted spaces}
The goal of this subsection is to prove a Proposition on the continuity of the operator $T_{\lambda,\eta}$ in a weighted functional space that we will define below. This continuity presents the key to the proof of the theorem on the existence of solutions. The proof is based on estimates of Green's function $K_{\lambda,\eta}$ with weights. We will start by defining the weights as well as the functional spaces on which we  work, then we give two lemmas on which is based the proof of the proposition on continuity.\\
Let $\eta_0, \lambda_0>0$ small enough and let $\eta\in[0,\eta_0]$, $\lambda\in\mathbb{C}$ such that $|\lambda|\leqslant \lambda_0$. We let $$\langle v \rangle = \sqrt{1+|v|^2}, \ v\in \mathbb{R} $$
and define the weights $p_i^{\lambda,\eta}$ for $i=1, 2$ by:
\begin{description}
\item [For $\eta>0$:]
\begin{equation}\label{p^eta} 
p_2^{\lambda,\eta}(v):= \left\{\begin{array}{l} \eta^{\frac{\gamma}{3}}|\mathfrak{a}_{\lambda,\eta}(v)|  ,   \quad v \geqslant s_0\eta^{-\frac{1}{3}} , \\
 \langle v \rangle^{-\gamma},\qquad  \quad |v| \leqslant s_0\eta^{-\frac{1}{3}} , \\
\eta^{\frac{\gamma}{3}}|\mathfrak{b}_{\lambda,\eta}(v)|  , \quad v \leqslant -s_0\eta^{-\frac{1}{3}}  .
\end{array}\right. 
\end{equation}
\item[For $\eta=0$:] 
\begin{equation}\label{p^0}
p_2^0(v):= \langle v \rangle^{-\gamma} , \quad\forall v\in \mathbb{R}
\end{equation}
\end{description}
and  
$$ p_1^{\lambda,\eta}(v):= \frac{p_2^{\lambda,\eta}(v)}{\langle v\rangle^{2+\delta}} , \quad\forall v\in \mathbb{R} , \forall \eta\in [0,\eta_0] ,$$ 
for any $\delta \in (0,2)$. \\
Note that $p_2^{\lambda,\eta}$ belongs to $L^2(\mathbb{R})$ since $2\gamma >1$.
We define the Banach space $E^\eta_i$, as the completion of $C^\infty_c(\mathbb{R},\mathbb{C})$ for the norm: $\ds \| \phi \|_{E^\eta_i} := \left\| \frac{\phi}{p_i^{\lambda,\eta}}\right\|_{L^\infty}$
$$ E^\eta_i := \overline{ \{ \phi \in C^\infty_c(\mathbb{R},\mathbb{C}) /\  \| \phi \|_{E^\eta_i} < +\infty\} }.$$ 
We have the embeddings
$$ \|\cdot\|_{E^{\eta^*}_i} \leqslant \|\cdot\|_{E_i^\eta} \ \mbox{for } \eta \leqslant \eta^*.$$
\begin{lemma}\label{continuite avec sigma}
Let $\delta\in(0,2)$. Then, there exists a constant $C>0$ such that:
\begin{equation}\label{continuite de T_0}
\int_\mathbb{R} |K_0(v,w)|\langle w\rangle^{-\gamma-\delta-2}\frac{\mathrm{d}w}{\langle v\rangle^{-\gamma}} \leqslant C, \quad \forall v \in \mathbb{R}.
\end{equation}
thus $T_0$ in $\mathcal{L}(E^0_1,E^0_2)$ is continuous.
\end{lemma}

\noindent \bp Thanks to the second point of the Remark \ref{remarques sur T_eta}, since the weights are symmetric, we establish the inequality \eqref{continuite de T_0} on $\RR^+$. Let $\delta\in(0,2)$ and let $v_0>0$ big enough.\\  
\noindent \textbf{Step 1: $v\geqslant v_0$.} We have
\begin{align*}
\int_\mathbb{R} |K_0(v,w)|\langle w\rangle^{-\gamma-\delta-2}\frac{\mathrm{d}w}{\langle v\rangle^{-\gamma}} &\leqslant \langle v\rangle^{\gamma} \bigg[ \psi_1(v)\bigg(\int_{-\infty}^{-v_0} \frac{\psi_2(w)\mathrm{d}w}{|w|^{\gamma+\delta+2}}+ \int_{-v_0}^{v_0}  \frac{\psi_2(w) \mathrm{d}w}{\langle w\rangle^{\gamma+\delta+2}} + \int_{v_0}^v  \frac{\psi_2(w) \mathrm{d}w}{w^{\gamma+\delta+2}}\bigg) \\
&+\psi_2(v)\int_v^\infty \frac{\psi_1(w) \mathrm{d}w}{\langle w\rangle^{\gamma+\delta+2}}\bigg] .
\end{align*}
For $v\geqslant v_0$: $\psi_1(v)\lesssim v^{-\gamma}$ and $\psi_2(v) \lesssim v^{\gamma+1}$,  on $[-v_0,v_0]$: $\psi_1$ and $\psi_2$ are bounded by a positive constant, and for $v\leqslant -v_0$ we have: $\psi_1(v)\lesssim |v|^{\gamma+1}$ and $\psi_2(v)\lesssim |v|^{-\gamma}$ since $\psi_1(-v)=\psi_2(v)$. Then,
\begin{align*}
\int_\mathbb{R} |K_0(v,w)|\langle w\rangle^{-\gamma-\delta-2}\frac{\mathrm{d}w}{\langle v\rangle^{-\gamma}}  &\lesssim \int_{-\infty}^{-v_0} \frac{\mathrm{d}w}{|w|^{2\gamma+\delta+2}}+ 1 + \int_{v_0}^v \frac{\mathrm{d}w}{|w|^{\delta+1}} + v^{2\gamma+1}\int_v^\infty \frac{\mathrm{d}w}{w^{2\gamma+\delta+2}} \\
&\lesssim 1 + v^{-\delta} \lesssim 1 .
\end{align*}
\noindent \textbf{Step 2: $v\in [0,v_0]$.} In this case, we cut the integral as follows:
\begin{align*}
\int_\mathbb{R} |K_0(v,w)|\langle w\rangle^{-\gamma-\delta-2}\frac{\mathrm{d}w}{\langle v\rangle^{-\gamma}} &\leqslant \langle v\rangle^{\gamma} \bigg[ \psi_1(v)\bigg(\int_{-\infty}^{-v_0} \frac{\psi_2(w)\mathrm{d}w}{|w|^{\gamma+\delta+2}} + \int_{-v_0}^v \frac{\psi_2(w) \mathrm{d}w}{\langle w\rangle^{\gamma+\delta+2}} \bigg) \\
&+ \psi_2(v)\bigg(\int_v^{v_0} \frac{\psi_1(w) \mathrm{d}w}{\langle w\rangle^{\gamma+\delta+2}} + \int_{v_0}^\infty \frac{\psi_1(w)\mathrm{d}w}{w^{\gamma+\delta+2}} \bigg)\bigg] ,
\end{align*}
and as in the previous step: since $v \in[0,v_0]$ then, $\langle v\rangle^{\gamma}$, $\psi_1$ and $\psi_2$ are bounded by a positive constant. Also, we have for $w\leqslant -v_0$: $\psi_2(w)\lesssim |w|^{-\gamma}$, and for $w \geqslant v_0$: $\psi_1(w)\lesssim w^{-\gamma}$. Therefore,
\begin{align*}
\int_\mathbb{R} |K_0(v,w)|\langle w\rangle^{-\gamma-\delta-2}\frac{\mathrm{d}w}{\langle v\rangle^{-\gamma}} &\lesssim \int_{-\infty}^{-v_0} \frac{\mathrm{d}w}{|w|^{2\gamma+\delta+2}} + 1 + \int_{v_0}^\infty \frac{\mathrm{d}w}{w^{2\gamma+\delta+2}}  \lesssim 1.
\end{align*} 
\ep
\begin{lemma}\label{lem Kp_1< c p_2} Let $\eta_0, \lambda_0>0$ small enough. There exists a constant $C > 0$ independent of $v$, $\eta$ and $\lambda$ such that
\begin{equation}\label{Kp_1<p_2}
\int_\mathbb{R} |K_{\lambda,\eta}(v,w)|p_1^{\lambda,\eta}(w)\mathrm{d}w \leqslant C  p_2^{\lambda,\eta}(v)
\end{equation}
holds for all $v \in \mathbb{R}$, $\eta\in[0,\eta_0]$ and $|\lambda|\leqslant\lambda_0$.
\end{lemma}
\bp The case $\eta=0$ is treated in the previous Lemma so, let $\eta\neq 0$. We will proceed as in the previous Lemma and since $\bar K_{\bar\lambda,\eta}(-v,-w)=K_{\lambda,\eta}(v,w)$ and $\bar p^{\bar\lambda,\eta}_i(-v)=p^{\lambda,\eta}_i(v)$ for $i=1, 2$ then, we show the inequality for $v\in\RR^+$.  Let denote by
$$ E^{\lambda,\eta}(v):=\int_\mathbb{R} |K_{\lambda,\eta}(v,w)|\frac{p_1^{\lambda,\eta}(w)}{p_2^{\lambda,\eta}(v)}\mathrm{d}w .$$
\noindent \textbf{Step 1: $v\geqslant s_0\eta^{-\frac{1}{3}}$.} We will cut the integral into four parts, according to the behaviors of the $\psi^{\lambda,\eta}_i$ and the definition of the weight, as follows:
\begin{align*}
E^{\lambda,\eta}(v) &\leqslant \frac{|\psi^{\lambda,\eta}_1(v)|}{p_2^{\lambda,\eta}(v)}\bigg[ \int_{-\infty}^{-s_0\eta^{-\frac{1}{3}}} |\psi^{\lambda,\eta}_2| p_1^{\lambda,\eta}(w)\mathrm{d}w+ \int_{-s_0\eta^{-\frac{1}{3}}}^{s_0\eta^{-\frac{1}{3}}} |\psi^{\lambda,\eta}_2|p_1^{\lambda,\eta}(w) \mathrm{d}w \\ 
&+ \int_{s_0\eta^{-\frac{1}{3}}}^v  |\psi^{\lambda,\eta}_2|p_1^{\lambda,\eta}(w) \mathrm{d}w\bigg] + \frac{|\psi^{\lambda,\eta}_2(v)|}{p_2^{\lambda,\eta}(v)} \int_v^\infty  |\psi^{\lambda,\eta}_1|p_1^{\lambda,\eta}(w) \mathrm{d}w .
\end{align*}
Thanks to the inequalities \eqref{psi_1^eta} and \eqref{psi_2^eta} of Lemma \ref{estimation des psi^eta_i}, we have for $v \geqslant s_0\eta^{-\frac{1}{3}}$,
$$ \frac{|\psi^{\lambda,\eta}_1(v)|}{p_2^{\lambda,\eta}(v)} \lesssim 1 , \ |\psi^{\lambda,\eta}_2(v)| \lesssim \eta^{-\frac{\gamma+1}{3}}|\mathfrak{b}_{\lambda,\eta}(v)| \  \mbox{ and }  \ \frac{|\psi^{\lambda,\eta}_2(v)|}{p_2^{\lambda,\eta}(v)} \lesssim \eta^{-\frac{2\gamma+1}{3}}\frac{|\mathfrak{b}_{\lambda,\eta}(v)|}{|\mathfrak{a}_{\lambda,\eta}(v)|} .$$
For $|w|\leqslant s_0\eta^{-\frac{1}{3}}$, $|\psi^{\lambda,\eta}_2(w)|\lesssim \psi_2(w)\lesssim \langle w\rangle^{\gamma+1}$ and finally for $w\leqslant -s_0\eta^{-\frac{1}{3}}$, $|\psi^{\lambda,\eta}_2(w)|\lesssim p_2^{\lambda,\eta}(w)$, with $p_2^{\lambda,\eta}(w)= \eta^{\frac{\gamma}{3}}|\mathfrak{a}_{\lambda,\eta}(w)|$ for this range of velocity.  Therefore, 
\begin{align*}
E^{\lambda,\eta}(v) &\lesssim  \int_{-\infty}^{-s_0\eta^{-\frac{1}{3}}} \eta^{\frac{2\gamma}{3}}|\mathfrak{b}_{\lambda,\eta}(w)|^2\frac{\mathrm{d}w}{\langle w\rangle^{2+\delta}}+ \int_{-s_0\eta^{-\frac{1}{3}}}^{s_0\eta^{-\frac{1}{3}}} \frac{\mathrm{d}w}{\langle w\rangle^{1+\delta}} \\ 
&+ \int_{s_0\eta^{-\frac{1}{3}}}^v  \eta^{-\frac{1}{3}} \frac{|\mathfrak{a}_{\lambda,\eta}(w)\mathfrak{b}_{\lambda,\eta}(w)|}{\langle w\rangle^{2+\delta}} \mathrm{d}w + \eta^{-\frac{1}{3}}\bigg|\frac{\mathfrak{b}_{\lambda,\eta}(v)}{\mathfrak{a}_{\lambda,\eta}(v)}\bigg| \int_v^\infty \frac{|\mathfrak{a}_{\lambda,\eta}(w)|^2}{\langle w\rangle^{2+\delta}} \mathrm{d}w \\
&= \sum_{i=1}^4 I^{\lambda,\eta}_i .
\end{align*}
In order to estimate $I^{\lambda,\eta}_1$, $I^{\lambda,\eta}_3$ and $I^{\lambda,\eta}_4$, we will perform the changes of variables $w=\eta^{-\frac{1}{3}}t$ and $v=\eta^{-\frac{1}{3}}s$. \\
$\bullet$ \textbf{Estimation of $I^{\lambda,\eta}_1$:} 
$$
I^{\lambda,\eta}_1:=\eta^{\frac{2\gamma}{3}} \int_{-\infty}^{-s_0\eta^{-\frac{1}{3}}} \frac{|\mathfrak{b}_{\lambda,\eta}(w)|^2}{\langle w\rangle^{2+\delta}}\mathrm{d}w =\eta^{\frac{2\gamma-1}{3}}  \int_{-\infty}^{-s_0} \frac{|\mathfrak{b}_\lambda(t)|^2}{\langle \eta^{-\frac{1}{3}}t \rangle^{2+\delta}}\mathrm{d}t \lesssim \eta^{\frac{2\gamma+\delta+1}{3}}  \int_{-\infty}^{-s_0} \frac{|\mathfrak{b}_\lambda(t)|^2}{ |t|^{2+\delta}}\mathrm{d}t \lesssim 1 .
$$
$\bullet$ \textbf{Estimation of $I^{\lambda,\eta}_2$:}
$$
I^{\lambda,\eta}_2:=\int_{-s_0\eta^{-\frac{1}{3}}}^{s_0\eta^{-\frac{1}{3}}} \frac{\mathrm{d}w}{\langle w\rangle^{1+\delta}} \leqslant \int_\RR \frac{\mathrm{d}w}{\langle w\rangle^{1+\delta}} \lesssim 1 .
$$
$\bullet$ \textbf{Estimation of $I^{\lambda,\eta}_3(v)$:} since $|\mathfrak{a}_\lambda(t)\mathfrak{b}_\lambda(t)|\lesssim t^{-\frac{1}{2}}$ for $|\lambda|\leqslant \lambda_0$ and $t\geqslant s_0$, then
$$
I^{\lambda,\eta}_3(v):=\eta^{-\frac{1}{3}} \int_{s_0\eta^{-\frac{1}{3}}}^v   \frac{|\mathfrak{a}_{\lambda,\eta}(w)\mathfrak{b}_{\lambda,\eta}(w)|}{\langle w\rangle^{2+\delta}} \mathrm{d}w= \eta^{-\frac{2}{3}} \int_{s_0}^s   \frac{|\mathfrak{a}_\lambda(t)\mathfrak{b}_\lambda(t)|}{\langle  \eta^{-\frac{1}{3}}t\rangle^{2+\delta}} \mathrm{d}t \lesssim \eta^{\frac{\delta}{3}} \int_{s_0}^s   \frac{\mathrm{d}t}{t^{\frac{5}{2}+\delta}} \lesssim 1  .
$$
$\bullet$ \textbf{Estimation of $I^{\lambda,\eta}_4(v)$:} we have
$$
 I^{\lambda,\eta}_4(v):=\eta^{-\frac{1}{3}}\bigg|\frac{\mathfrak{b}_{\lambda,\eta}(v)}{\mathfrak{a}_{\lambda,\eta}(v)}\bigg| \int_v^\infty \frac{|\mathfrak{a}_{\lambda,\eta}(w)|^2}{\langle w\rangle^{2+\delta}} \mathrm{d}w = \eta^{-\frac{2}{3}}\bigg|\frac{\mathfrak{b}_\lambda(s)}{\mathfrak{a}_\lambda(s)}\bigg| \int_s^\infty \frac{|\mathfrak{a}_\lambda(t)|^2}{\langle \eta^{-\frac{1}{3}}t\rangle^{2+\delta}} \mathrm{d}s \lesssim \frac{\eta^{\frac{\delta}{3}}}{s^\frac{1}{2}} \int_{s_0}^s   \frac{\mathrm{d}t}{t^{2+\delta}} \lesssim 1 
$$ 
since $\frac{|\mathfrak{a}_\lambda(t)|}{|\mathfrak{a}_\lambda(s)|}\leqslant 1$ and $|\mathfrak{a}_\lambda(s)\mathfrak{b}_\lambda(s)|\lesssim t^{-\frac{1}{2}}$ for  $t\geqslant s \geqslant s_0$ and $|\lambda|\leqslant \lambda_0$ (by the inequalities of \eqref{estimations a_lambda b_lambda}). \\

\noindent \textbf{Step 2: $v\in[0,s_0\eta^{-\frac{1}{3}}]$.} We will proceed exactly as in the previous step by cutting the integral this time as follows
\begin{align*}
E^{\lambda,\eta}(v) &\leqslant \frac{|\psi^{\lambda,\eta}_1(v)|}{p_2^{\lambda,\eta}(v)}\bigg[ \int_{-\infty}^{-s_0\eta^{-\frac{1}{3}}} |\psi^{\lambda,\eta}_2| p_1^{\lambda,\eta}(w)\mathrm{d}w+ \int_{-s_0\eta^{-\frac{1}{3}}}^v |\psi^{\lambda,\eta}_2|p_1^{\lambda,\eta}(w) \mathrm{d}w \\ 
&+ \int_v^{s_0\eta^{-\frac{1}{3}}}  |\psi^{\lambda,\eta}_2|p_1^{\lambda,\eta}(w) \mathrm{d}w\bigg] + \frac{|\psi^{\lambda,\eta}_2(v)|}{p_2^{\lambda,\eta}(v)} \int_{s_0\eta^{-\frac{1}{3}}}^\infty  |\psi^{\lambda,\eta}_1|p_1^{\lambda,\eta}(w) \mathrm{d}w  .
\end{align*}
We have in this case:  $p_2^{\lambda,\eta}(v) = \langle v \rangle^{-\gamma}$ and $|\psi^{\lambda,\eta}_i(v)|\lesssim \psi_i(v)$ for $i=1,2$ by the inequalities \eqref{psi_1^eta} and \eqref{psi_2^eta}, with $\psi_1(v)\lesssim p_2^{\lambda,\eta}(v)$ and $\psi_2(v)\lesssim \langle v\rangle^{\gamma+1}$ by \eqref{comportement de psi_1} and \eqref{psi_2}. Then,   
\begin{align*}
E^{\lambda,\eta}(v) &\lesssim  \int_{-\infty}^{-s_0\eta^{-\frac{1}{3}}} \eta^{\frac{2\gamma}{3}}|\mathfrak{b}_{\lambda,\eta}(w)|^2\frac{\mathrm{d}w}{\langle w\rangle^{2+\delta}}+ \int_{-s_0\eta^{-\frac{1}{3}}}^v \frac{\mathrm{d}w}{\langle w\rangle^{1+\delta}} \\ 
&+ \langle v\rangle^{2\gamma+1} \int_v^{s_0\eta^{-\frac{1}{3}}}  \frac{\mathrm{d}w}{\langle w\rangle^{2\gamma+\delta+2}}  + \eta^{\frac{\gamma}{3}}\langle v\rangle^{2\gamma+1} \int_{s_0\eta^{-\frac{1}{3}}} ^\infty \frac{|\mathfrak{a}_{\lambda,\eta}(w)|}{\langle w\rangle^{\gamma+\delta+2}} \mathrm{d}w  .
\end{align*}
The first two integrals are bounded by $I^{\lambda,\eta}_1+I^{\lambda,\eta}_2$ which is uniformly bounded with respect to $v$ by step 1. For the last two terms, we write
$$
\langle v\rangle^{2\gamma+1} \int_v^{s_0\eta^{-\frac{1}{3}}}  \frac{\mathrm{d}w}{\langle w\rangle^{2\gamma+\delta+2}} = \int_v^{s_0\eta^{-\frac{1}{3}}}  \frac{\langle v\rangle^{2\gamma+1}}{\langle w\rangle^{2\gamma+1}}\frac{\mathrm{d}w}{\langle w\rangle^{1+\delta}} \leqslant  \int_\RR  \frac{\mathrm{d}w}{\langle w\rangle^{1+\delta}} \lesssim 1 
$$
and 
$$
\eta^{\frac{\gamma}{3}}\langle v\rangle^{2\gamma+1} \int_{s_0\eta^{-\frac{1}{3}}} ^\infty \frac{|\mathfrak{a}_{\lambda,\eta}(w)|}{\langle w\rangle^{\gamma+\delta+2}} \mathrm{d}w = \eta^{\frac{\gamma}{3}}\langle v\rangle^{\gamma} \int_{s_0\eta^{-\frac{1}{3}}}^\infty \frac{\langle v\rangle^{\gamma+1}}{\langle w\rangle^{\gamma+1}}\frac{|\mathfrak{a}_{\lambda,\eta}(w)|}{\langle w\rangle^{1+\delta}} \mathrm{d}w \lesssim  \int_\RR  \frac{\mathrm{d}w}{\langle w\rangle^{1+\delta}} \lesssim 1 
$$
since $v\leqslant s_0\eta^{-\frac{1}{3}}$ and $|\mathfrak{a}_{\lambda,\eta}(w)| \leqslant 1$ for $w \geqslant s_0\eta^{-\frac{1}{3}}$. 
\ep
\begin{proposition}[Continuity of $T_{\lambda,\eta}$]\label{continuite de T_eta}
Let $\eta_0, \lambda_0>0$ small enough. Let $\eta \in [0,\eta_0]$ and $\lambda \in \mathbb{C}$ such that $|\lambda|\leqslant\lambda_0$. Then, the map $T_{\lambda,\eta} : E^\eta_1 \longrightarrow E^\eta_2$ is linear continous. \\
Moreover, there exists $C > 0$ independent of $\eta$ and $\lambda$ such that
\begin{equation}\label{|T(g)|<|g|}
\| T_{\lambda,\eta}(g) \|_{E^\eta_2} \leqslant C \|g\|_{E^\eta_1} , \quad \forall g \in E^\eta_1 .
\end{equation}
\end{proposition}
\bp The proof of this Proposition is a direct consequence of Lemma \ref{lem Kp_1< c p_2}.
\ep
\section{Existence of the eigen-solution}
\subsection{Existence of solutions for the penalized equation}
In this subsection, we use the ``right inverse" operator $T_{\lambda,\eta}$ to rewrite once again the penalized equation \eqref{eq penalisee2} as a fixed point problem for the identity plus a compact map. Then, the Fredholm Alternative will allow us to apply the Implicit Function Theorem in order to get existence of solutions for this new problem, thus solutions for the equation \eqref{eq penalisee2}.\\
Define $F: \{\lambda\in\CC;|\lambda|\leqslant\lambda_0\}\times[0,\eta_0]\times C_b(\RR,\CC)\longrightarrow C_b(\RR,\CC)$ by 
$$F(\lambda,\eta,h):= h - \mathcal{T}_{\lambda,\eta}(h) ,$$ 
with
$$\mathcal{T}_{\lambda,\eta}(h):=\frac{1}{p_2^{\lambda,\eta}}T_{\lambda,\eta}\big[Vp_2^{\lambda,\eta}h-\langle p_2^{\lambda,\eta}h-M,\Phi \rangle \Phi \big] .$$
Note that finding a solution $h(\lambda,\eta)$ solution to $F\big(\lambda,\eta,h(\lambda,\eta)\big)= 0$ gives a solution to the penalized equation by taking $M_{\lambda,\eta}= h(\lambda,\eta) p_2^{\lambda,\eta}$. \\

\ni The function $\Phi$ is a function satisfying the following assumptions:
\begin{enumerate}
\item For all $v$ in $\RR$, $\Phi(v)=\Phi(-v)>0$.
\item For all $\eps>0$, there exists $g^\eps_1 \in C_c(\RR)$ such that $\| \Phi/p^{\lambda,\eta}_1 - g^\eps_1 \|_\infty < \eps$ with $\mathrm{supp}(g^\eps_1)$ independent of $\lambda$ and $\eta$. 
\item Even if it means multiplying $\Phi$ by a constant, we can take it such that $\langle \Phi, M \rangle = 1$.
\end{enumerate}
Any continuous function with compact support and satisfying 1. and 3. is suitable. The function $\Phi:=\Phi_{\lambda,\eta}=c_{\lambda,\eta} \langle v\rangle^{-\sigma} p^{\lambda,\eta}_1$ satisfies all the previous assumptions, where $\sigma>0$ and $c_{\lambda,\eta}$ is a constant of ``normalization" such that $\langle \Phi_{\lambda,\eta},M\rangle = 1$. \\

\begin{remark}\label{Phi et T_lambda,eta}
\item \begin{enumerate}
\item For the following, we fix a continuous function $\Phi$, with compact support included in $[-R,R]$ with $R > 0$, and satisfying assumptions 1. and 3. above.
\item Note that $\mathcal{T}_{\lambda,0}$ does not depend on $\lambda$ since $T_{\lambda,0}$ does not. Let's denote it by $\mathcal{T}_0$:
$$\mathcal{T}_0(h):=\mathcal{T}_{\lambda,0}(h)=\frac{1}{p_2^0}T_{0}\big[Vp_2^0 h-\langle p_2^0 h-M,\Phi\rangle \Phi\big] $$
\item The map $\mathcal{T}_{\lambda,\eta}$ is affine with respect to $h$. We denote by $\mathcal{T}^l_{\lambda,\eta}$ its linear part:
$$ \mathcal{T}^l_{\lambda,\eta}(h):=\frac{1}{p_2^{\lambda,\eta}}T_{\lambda,\eta}\big[V p_2^{\lambda,\eta}h-\langle p_2^{\lambda,\eta}h,\Phi\rangle \Phi\big] .$$
\end{enumerate}
\end{remark}
\begin{lemma}[Continuity, differentiability and compactness of $\mathcal{T}_{\lambda,\eta}$]\label{conti. de T_eta en eta}
 Let $\eta_0, \lambda_0>0$ small enough. Let $\eta \in [0,\eta_0]$ and $\lambda \in \mathbb{C}$ such that $|\lambda|\leqslant\lambda_0$. Then, the map $\mathcal{T}_{\lambda,\eta}:C_b(\mathbb{R},\mathbb{C})\longrightarrow C_b(\mathbb{R},\mathbb{C})$ is 
\begin{enumerate}
\item continuous with respect to $\lambda$ and $\eta$, moreover
\begin{equation}\label{T_eta - T_eta' -->0}
\|\mathcal{T}_{\lambda,\eta}-\mathcal{T}_{\lambda,\eta'}\|_{\mathcal{L}(C_b)} \underset{\eta\rightarrow \eta'}{\longrightarrow} 0, \quad \forall \eta' \in [0,\eta_0]
\end{equation}
and 
\begin{equation}\label{T_lambda - T_lambda' --> 0}
\|\mathcal{T}_{\lambda,\eta}-\mathcal{T}_{\lambda',\eta}\|_{\mathcal{L}(C_b)} \underset{\lambda \rightarrow \lambda'}{\longrightarrow} 0 , \quad \forall \lambda' \in \CC ; |\lambda'| \leqslant \lambda_0.
\end{equation}
\item differentiable in $C_b(\mathbb{R},\mathbb{C})$ and its differential is
\begin{equation}\label{diff T_eta}
 \frac{\partial\mathcal{T}_{\lambda,\eta}}{\partial h}=\mathcal{T}^l_{\lambda,\eta}=\frac{1}{p_2^{\lambda,\eta}}T_{\lambda,\eta}\big[Vp_2^{\lambda,\eta}\cdot-\langle p_2^{\lambda,\eta}\cdot,\Phi\rangle \Phi\big].
\end{equation}
\item The map $\mathcal{T}^l_0$ is compact.
\end{enumerate}
where $\mathcal{L}(C_b)$ is the space of linear continuous operators from $C_b(\RR,\CC)$ to itself.
\end{lemma}
\bpl \ref{conti. de T_eta en eta}. \\
1. Continuity of $\mathcal{T}_{\lambda,\eta}$ with respect to $\lambda$ and $\eta$.  It is sufficient to prove this continuity for the map $\mathcal{T}_{\lambda,\eta}$ composed with the characteristic function $\chi_{[-R,R]}$, since $\mathcal{T}_{\lambda,\eta}^l$ can be written as 
$$ \mathcal{T}_{\lambda,\eta}^l(h)=\frac{1}{p^{\lambda,\eta}_2}T_{\lambda,\eta}\big[ g_1 p^{\lambda,\eta}_1 h - \langle p_2^{\lambda,\eta}h,\Phi\rangle g_2 p^{\lambda,\eta}_1\big] , $$ 
with $g_1:=Vp^{\lambda,\eta}_2/p^{\lambda,\eta}_1$ and $g_2:=\Phi/p^{\lambda,\eta}_1$ belong to $C_0(\RR)$, the set of continuous functions converging to zero at the infinity. Indeed, let us denote by $\mathcal{T}_{\lambda,\eta}^R:=\mathcal{T}_{\lambda,\eta}^l\circ\chi_{[-R,R]}$. Let $h\in C_b(\RR,\CC)$ and let $\varepsilon>0$. Then, there exists $g^\varepsilon_1, g^\varepsilon_2 \in C^\infty_c(\RR)$ such that $\|g^\varepsilon_i-g_i\|_\infty \leqslant \varepsilon/(2C)$ for $i=1,2$, where $C$ is the constant of Lemma \ref{lem Kp_1< c p_2}. Let $R_\varepsilon>0$ such that $\mathrm{supp}(g^\varepsilon_1)\cup\mathrm{supp}(g^\varepsilon_2)\subset [-R_\varepsilon,R_\varepsilon]$. We have thanks to the Lemma \ref{lem Kp_1< c p_2} and since $p_2^{\lambda,\eta}(v)\leqslant p_2^0(v)$ with $\langle p_2^0,\Phi\rangle=1$:
\begin{align*}
\| \mathcal{T}_{\lambda,\eta}^l(h)-\mathcal{T}_{\lambda,\eta}^{R_\varepsilon}(h)\|_\infty &\leqslant C \big[\|g^\varepsilon_1-g_1\|_\infty\|h\|_\infty+\|g^\varepsilon_2-g_2\|_\infty|\langle p_2^{\lambda,\eta}h,\Phi\rangle|\big] \\
&\leqslant C \big[\|g^\varepsilon_1-g_1\|_\infty+\|g^\varepsilon_2-g_2\|_\infty|\langle p_2^0,\Phi\rangle|\big]\|h\|_\infty .
\end{align*}
Hence, 
\begin{equation}\label{T_eta - T_eta^R}
\| \mathcal{T}_{\lambda,\eta}^l(h)-\mathcal{T}_{\lambda,\eta}^{R_\varepsilon}(h)\|_\infty \leqslant \varepsilon \|h\|_\infty , \quad \forall h \in C_b(\RR,\CC) .
\end{equation}
\begin{remark} 
\item \begin{enumerate}
\item Note that $\mathrm{supp}(g^\eps_1)$ does not depend on $h$, $\lambda$ and $\eta$, since $g_1$ does not depend on the latter three: $g_1 = \langle v \rangle^{2+\delta}V = \frac{\gamma(\gamma+2)}{\langle v \rangle^{2-\delta}} \in C_0 (\RR,\CC)$ for $\delta\in(0,2)$. Similarly for $g_2$, by assumption. Moreover, since $\Phi$ has compact support then, for $s_0$ large enough and $\eta$ small enough, $\mathrm{supp}(\Phi) \subset [-s_0\eta^{-\frac{1}{3}},s_0\eta^{-\frac{1}{3}}]$ with $p^{\lambda,\eta}_1(v) = p^0_1(v)$ on this last interval. Therefore, $R_\varepsilon$ is independent of $h$, $\lambda$ and $\eta$. 
\item Since $\Phi$ has compact support then,  $|\langle p^{\lambda,\eta}_2,\Phi\rangle|=|\langle p^0_2,\Phi\rangle| = 1$ for $\eta$ small enough. Therefore,  
$$ \big\|\mathcal{T}_{\lambda,\eta}(\langle p^{\lambda,\eta}_2 h,\Phi\rangle \Phi) \big\|_\infty \lesssim \|h\|_\infty.$$
\end{enumerate} 
\end{remark}
Let us now show the continuity of $\mathcal{T}_{\lambda,\eta}^R$. Let $|v|\leqslant R$ and let $\eta_0$ small enough such that $\eta_0^{-\frac{1}{3}} > R$.  Then, $p^{\lambda,\eta}_2(v)=p^0_2(v)$ for all $\eta \in [0,\eta_0]$ and for all $|v| \leqslant R$, and we have: 
$$
\mathcal{T}_{\lambda,\eta}^R(h)(v)-\mathcal{T}_{\lambda,\eta'}^R(h)(v)=\int_{-R}^R \big[K_{\lambda,\eta}(v,w)-K_{\lambda,\eta'}(v,w)\big]\big[g_1 h(w)-\langle p^0_2h,\Phi\rangle g_2(w)\big]\frac{p^0_1(w)}{p^0_2(v)}\mathrm{d}w .
$$
Therefore,  
$$
|\mathcal{T}_{\lambda,\eta}^R(h)(v)-\mathcal{T}_{\lambda,\eta'}^R(h)(v)| \leqslant \int_{-R}^R  \big| K_{\lambda,\eta}(v,w)-K_{\lambda,\eta'}(v,w) \big| \frac{p^0_1(w)}{p^0_2(v)}\mathrm{d}w\big(\|g_1\|_\infty+\|g_2\|_\infty\big) \|h\|_\infty . $$
Similarly for $\lambda$, 
$$ |\mathcal{T}_{\lambda,\eta}^R(h)(v)-\mathcal{T}_{\lambda',\eta}^R(h)(v)| \leqslant \int_{-R}^R  \big| K_{\lambda,\eta}(v,w)-K_{\lambda',\eta}(v,w) \big| \frac{p^0_1(w)}{p^0_2(v)}\mathrm{d}w\big(\|g_1\|_\infty+\|g_2\|_\infty\big) \|h\|_\infty . $$
We conclude with Lebesgue's theorem thanks to the continuity of $\| K_{\lambda,\eta}\|_{L^\infty([-R,R]\times[-R,R])}$ with respect to $\lambda$ and $\eta$,  and since $|K_{\lambda,\eta}|$ is dominated by $K_0$ on $[-R,R]\times[-R,R] $ thanks to \eqref{K_eta < K_0}.
Hence the limits \eqref{T_eta - T_eta' -->0} and \eqref{T_lambda - T_lambda' --> 0} hold. 
\\
2. The second point is immediate since $\mathcal{T}_{\lambda,\eta}$ is an affine map with respect to $h$.\\
3. Let us first show that $\mathcal{T}_0^R$ is compact by  Ascoli-Arz\'{e}la theorem. Define $B_R(0,1):=\{h\in C([-R,R],\CC); \|h\|_\infty\leqslant1\}$ and introduce $\mathcal{F}:=\mathcal{T}_0^R(B_R(0,1)) $. The set $\mathcal{F}$ is bounded because $\mathcal{T}_0^R$ is bounded. The set $\mathcal{F}$ is equicontinuous since  for $|v_1-v_2|\leqslant \varepsilon/C_R$ we have: $|\mathcal{T}_0^R(h)(v_1)-\mathcal{T}_0^R(h)(v_2)|\leqslant \varepsilon$, $\forall h\in B_R(0,1)$, where and $C_R$ is a constant that depends only on $R$. Indeed,  
\begin{align*}
\mathcal{T}_0^R(h)(v_1)-\mathcal{T}_0^R(h)(v_2) &= \int_{|w|\leqslant R} \bigg(\frac{K_0(v_1,w)}{p^0_2(v_1)}-\frac{K_0(v_2,w)}{p^0_2(v_2)}\bigg)\big[g_1 h-\langle p^0_2h,\Phi\rangle g_2\big](w)p^0_1(w)\mathrm{d}w \\
&= \int_{|w|\leqslant R} \frac{p^0_2(v_2)-p^0_2(v_1)}{p^0_2(v_2)}K_0(v_1,w)\frac{p^0_1(w)}{p^0_2(v_1)}\big[g_1 h-\langle p^0_2h,\Phi\rangle g_2\big](w)\mathrm{d}w \\
&+ \int_{|w|\leqslant R} \big[K_0(v_1,w)-K_0(v_2,w)\big]\frac{p^0_1(w)}{p^0_2(v_2)}\big[g_1 h-\langle p^0_2h,\Phi\rangle g_2\big](w)\mathrm{d}w \\
&=: I_1+I_2 .
\end{align*}
Since $p^0_2$ is Lipschitz on the compact $[R,R]$ then,
$$ |I_1| \leqslant \tilde C_R |v_1-v_2| \big(\|g_1\|_\infty+\|g_2\|_\infty\big)\|h\|_\infty \leqslant C_R' |v_1-v_2| . $$
For $I_2$ we write:
\begin{align*}
I_2 &=  \int_{|w|\leqslant R}  \bigg[[\psi_1(v_1)-\psi_1(v_2)]\psi_2(w)\chi_{w<v_1}+\psi_1(w)[\psi_2(v_1)-\psi_2(v_2)]\chi_{w>v_1}\bigg]\frac{p^0_1(w)}{p^0_2(v_2)}\\
&\times\big[g_1 h-\langle p^0_2h,\Phi\rangle g_2\big](w)\mathrm{d}w + \int_{v_1}^{v_2} K_0(v_2,w)\frac{p^0_1(w)}{p^0_2(v_1)}\big[g_1 h-\langle p^0_2h,\Phi\rangle g_2\big](w)\mathrm{d}w .
\end{align*}
Thus, since $\psi_1$ and $\psi_2$ are Lipschitz on the compact $[R,R]$ then, $|I_2| \leqslant C_R'' |v_1-v_2|$.  Hence,
$$ \big|\mathcal{T}_0^R(h)(v_1)-\mathcal{T}_0^R(h)(v_2)\big| \leqslant C_R |v_1-v_2| .$$
The compactness of $\mathcal{T}_0^l$ follows from the compactness of $\mathcal{T}_0^R$ and the inequality \eqref{T_eta - T_eta^R}. Indeed, we have by \eqref{T_eta - T_eta^R}, $\|\mathcal{T}^l_0(h)-\mathcal{T}^{R_\eps}_0(h)\|_\infty \leqslant \eps,$ $\forall h\in C_b(\RR,\CC)$; $\|h\|_\infty \leqslant 1$, and since $\mathcal{T}^{R_\eps}_0$ is compact with $R_\eps$ being fixed and independent of $h$, $\lambda$ and $\eta$ then, there exists $N_\varepsilon \in \NN$, $\{h_i\}_{i=1}^{N_\varepsilon}\subset C_b(\RR,\CC)$ such that: $\mathcal{T}_0^{R_\eps}(h) \in \bigcup_{i=1}^{N_\varepsilon}B(h_i,\eps)$. Therefore, $\mathcal{T}_0^l(h)\in \bigcup_{i=1}^{N_\varepsilon}B(h_i,2\eps)$. Hence the compactness of $\mathcal{T}_0^l$ holds.
\epl
\begin{proposition}[Assumptions of the Implicit Function Theorem]\label{hypotheses de TFI}
\item \begin{enumerate}
\item The map $F(\lambda,\eta,\cdot)=Id-\mathcal{T}_{\lambda,\eta}$ is continuous in $C_b(\mathbb{R},\mathbb{C})$ uniformly with respect to $\lambda \mbox{ and } \eta$. Moreover, there exists $c>0$, independent of $\lambda$ and $\eta$ such that
$$ \|F(\lambda,\eta,h_1)-F(\lambda,\eta,h_2)\|_\infty \leqslant c \|h_1-h_2\|_\infty, \quad \forall h_1, h_2 \in C_b(\RR,\CC), \forall\eta, \forall |\lambda|\leqslant \lambda_0 .$$
\item $F$ is continuous with respect to $\lambda$ and $\eta$ and we have:
$$ \|F(\lambda,\eta,\cdot)-F(\lambda,\eta',\cdot)\|_{\mathcal{L}(C_b)} \underset{\eta \rightarrow \eta'}{\longrightarrow} 0 \ \mbox{ and }  \  \|F(\lambda,\eta,\cdot)-F(\lambda',\eta,\cdot)\|_{\mathcal{L}(C_b)} \underset{\lambda \rightarrow \lambda'}{\longrightarrow} 0 .$$
\item $F(\lambda,\eta,\cdot)$ is differentiable in $C_b(\mathbb{R},\mathbb{C})$, moreover:
 $$ \frac{\partial F}{\partial h}(\lambda,\eta,\cdot)=Id-\mathcal{T}_{\lambda,\eta}^l, \quad \forall |\lambda|\leqslant\lambda_0,\forall \eta\in[0,\eta_0] .$$
\item We have: $ F(0,0,\frac{M}{p^0_2})=0$ and $ \frac{\partial F}{\partial h}(0,0,\frac{M}{p^0_2})$ is invertible.
\end{enumerate}
\end{proposition}
\bpp \ref{hypotheses de TFI}.  \\
1. Let $h_1, h_2 \in C_b$ and let $\eta>0$ and $\lambda\in\mathbb{C}$ such that $|\lambda|\leqslant\lambda_0$. Then,
\begin{align*}
\|F(\lambda,\eta,h_1)-F(\lambda,\eta,h_2)\|_\infty &\leqslant \|(h_1-h_2)+\mathcal{T}_{\lambda,\eta}^l(h_1-h_2)\|_\infty \\
&\leqslant \big(1+C[\|g_1\|_\infty+\|g_2\|_\infty]\big)\|h_1-h_2\|_\infty\\
&\leqslant c\|h_1-h_2\|_\infty.
\end{align*}
2. The proof of this point is a direct consequence of the first point of Lemma \ref{conti. de T_eta en eta}.  \\
3. Follows from the second point of Lemma \ref{conti. de T_eta en eta}.  \\
4. We have for $(\lambda,\eta,h)=(0,0,M/p^0_2)$, $\ \langle p_2^0 h - M,\Phi\rangle = 0$. Then,
$$
F(0,0,M/p_2^0)=\frac{1}{p_2^0}\big(M-T_0[VM]\big) .
$$
Thus, multiplying this last equality by $p^0_2=M$ and applying the left inverse $L_0:=-\pa_v^2+\tilde W(v)$,  we obtain:
$$ L_0\big(M-T_0[VM]\big)=L_0(M)-VM = [-\pa_v^2+W(v)]M=0 . $$
Hence,  $F(0,0,\frac{M}{p^0_2})=0$ thanks to the injectivity of the left inverse $L_0$. \\
For the differential, we have $\frac{\partial F}{\partial h}(0,0,\frac{M}{p_2^0})=Id-\mathcal{T}_0^l$. By the Fredholm Alternative, this point is true if $\mathrm{Ker}(Id-\mathcal{T}_0^l)=\{0\}$. Let $h\in C_b(\mathbb{R},\mathbb{C})$ such that $h-\mathcal{T}_0^l(h)=0$.
By multiplying this last equation by $p_2^0$ and applying the operator $L_0$, we obtain
$$[-\partial_v^2+W(v)](p_2^0 h)=\langle p_2^0 h,\Phi\rangle\Phi . $$
Now, integrating the previous equation against $M$ and using the fact that $\langle \Phi,M\rangle=1$, we get 
$$ \langle p_2^0 h,\Phi\rangle = 0 .$$
Therefore, $p^0_2h$ is solution to $[-\partial_v^2+W(v)]f=0$. Then, there exists $c_1, c_2 \in \CC$ such that $p_2^0 h=c_1M+c_2Z$, which implies that $h=c_1\frac{M}{p_2^0}+c_2\frac{Z}{p_2^0} $. Since $h\in C_b$ and $\frac{Z}{p_2^0}\notin C_b$ then, $c_2=0$ and $h=c_1\frac{M}{p_2^0}$. Thus, $\langle p_2^0 h,\Phi\rangle =c_1=0$.  Hence, $h=0$. This completes the proof of the Proposition.
\epp
\begin{theorem}[Existence of solutions with constraint]\label{thm d'existence}
There is a unique function $M_{\lambda,\eta}$ in $E^\eta_2 \subset L^2(\RR,\CC)$ solution to
\begin{equation}\label{eq de M_eta + cont.}
[-\partial^2_v+ W(v)+ \mathrm{i} \eta v -\lambda\eta^{\frac{2}{3}}]M_{\lambda,\eta}(v)= b(\lambda,\eta)\Phi(v), \quad v\in\RR .
\end{equation} 
Moreover,
\begin{equation}\label{h_lambda,eta-h_0-->0}
\left\|\frac{M_{\lambda,\eta}}{p^{\lambda,\eta}_2}-\frac{M}{p^0_2}\right\|_{\infty} \underset{\eta \rightarrow 0}{\longrightarrow}0,
\end{equation}
where $\ds b(\lambda,\eta):=\langle N_{\lambda,\eta},\Phi\rangle$ with $N_{\lambda,\eta} := M_{\lambda,\eta}-M$.
\end{theorem}
\begin{remark}\label{rmq sur thm d'existence}
\item \begin{enumerate}
\item By construction, the solution $M_{\lambda,\eta}$ is symmetric and we have 
\begin{equation}\label{bar M_eta(-v)=M_eta(v)}
\bar{M_{\bar \lambda,\eta}}(-v) = M_{\lambda,\eta}(v) , \quad \forall v \in \RR.
\end{equation}
\item By introducing the function $c(\lambda,\eta)$, satisfying $c(\lambda,\eta)M_{\lambda,\eta}(0)=1$, we can always take $M_{\lambda,\eta}(0)=1$ in order to simplify the notations. Such a function $c(\lambda,\eta)$ exists since the solution $M_{\lambda,\eta}$ given by the theorem \ref{thm d'existence} does not vanish at $0$. Indeed, $M_{\lambda,\eta}(0)=0$ leads us to the following contradiction $$M(0)=1=\frac{M(0)}{p^0_2(0)}-\frac{M_{\lambda,\eta}(0)}{p^{\lambda,\eta}_2(0)}\leqslant\left\|\frac{M_{\lambda,\eta}}{p^{\lambda,\eta}_2}-\frac{M}{p^0_2}\right\|_{\infty} \underset{\eta \rightarrow 0}{\longrightarrow}0.$$
Moreover, $c(\lambda,\eta)$ is holomorphic in $\{|\lambda|\leqslant \lambda_0\}$ and continuous in $\eta\in[0,\eta_0]$. 
\end{enumerate}
\end{remark}
\bp By Proposition \ref{hypotheses de TFI}, $F$ satisfies the assumptions of the Implicit Function Theorem (IFT) around the point $(0,0,\frac{M}{p^0_2})$. Then, there exists $\lambda_0, \eta_0>0$ small enough, there exists a unique function $h:\{|\lambda|\leqslant\lambda_0\}\times\{|\eta|\leqslant\eta_0\}\longrightarrow C_b(\mathbb{R},\mathbb{C})$, continuous with respect to $\lambda$ and $\eta$ such that \begin{center}
$ F(\lambda,\eta,h(\lambda,\eta))=0$, for all $(\lambda,\eta) \in \{|\lambda|<\lambda_0\}\times\{|\eta |<\eta_0\}$.
\end{center}
Let denote $h_{\lambda,\eta}:=h(\lambda,\eta)$. Note that $h_{\lambda,0}$ does not depend on $\lambda$. The continuity of $h$ with respect to $\eta$ implies that 
\begin{equation}\label{h_eta-h_0-->0 dans L^infini}
\underset{\eta\rightarrow 0}{\lim}\|h_{\lambda,\eta}-h_{\lambda,0}\|_\infty=\underset{\eta\rightarrow 0}{\lim} \|h_{\lambda,\eta}-h_{0,0}\|_\infty = 0 .
\end{equation}
Finally, we take $M_{\lambda,\eta}:= p^{\lambda,\eta}_2 h_{\lambda,\eta}$ and the proof of the theorem is complete.
\ep
\subsection{Properties of the solution to the penalized equation}
\begin{corollary}[Properties of $M_{\lambda,\eta}$]\label{Propr de M_lambda,eta}
\item
\begin{enumerate}
\item There exists a constant $C$ such that, for all $v\in \mathbb{R}, |\lambda|\leqslant\lambda_0$ and $\eta\in[0,\eta_0]$\begin{equation}\label{M_eta < c M}
|M_{\lambda,\eta}(v)| \leqslant C M(v).
\end{equation}
\item For all $v\in \mathbb{R}$ and $|\lambda|\leqslant\lambda_0$ 
\begin{equation}\label{limite simple M_eta(v) = M(v)}
\underset{\eta \rightarrow 0}{\lim} M_{\lambda,\eta}(v) = M(v) .
\end{equation} 
Therefore,
\begin{equation}\label{cv dans L^2}
 \underset{\eta \rightarrow 0}{\lim} \int_\mathbb{R} M_{\lambda,\eta}(v)M(v)\mathrm{d}v = \int_\mathbb{R} M^2(v)\mathrm{d}v \quad \mbox{and } \ M_{\lambda,\eta}\underset{\eta \rightarrow 0}{\longrightarrow} M \mbox{ in } L^2(\RR).
\end{equation}
 \item  We have the following limit
 \begin{equation}\label{int eta^1:3 v M_eta M->0}
 \underset{\eta \rightarrow 0}{\lim} \int_\mathbb{R} \eta^{\frac{1}{3}}v M_{\lambda,\eta}(v)M(v)\mathrm{d}v = 0 .
\end{equation}  
 \end{enumerate}
\end{corollary}
\bp 1. We have
$$|M_{\lambda,\eta}(v)|\leqslant p_2^{\lambda,\eta}(v)\|h_{\lambda,\eta}\|_{\infty} \leqslant C p_2^0(v) = C M(v) $$
since $M_{\lambda,\eta}=p_2^{\lambda,\eta}h_{\lambda,\eta}$, with $\| h_{\lambda,\eta}\|_\infty \leqslant C$ uniformly with respect to $\lambda$ and $\eta$ thanks to \eqref{h_eta-h_0-->0 dans L^infini}, and since $p_2^{\lambda,\eta}(v)\leqslant p_2^0(v)$. For this last inequality, we have $p_2^{\lambda,\eta}(v)= p_2^0(v)$ for $v\in [-s_0\eta^{-\frac{1}{3}},s_0\eta^{-\frac{1}{3}}]$, and for $|v|\geqslant s_0\eta^{-\frac{1}{3}}$ we have: $ p_2^{\lambda,\eta}(v) / p_2^0(v) \leqslant (\eta^{-\frac{1}{3}}v)^{-\gamma}e^{-\frac{\sqrt{2}}{3}(\eta^{-\frac{1}{3}}v)^{\frac{3}{2}}} \leqslant 1$ since the function $t\mapsto t^{-\gamma}e^{-\frac{\sqrt{2}}{3}t^{\frac{3}{2}}}$ is decreasing for $t\in [s_0,+\infty)$ since $\gamma>0$. Hence, the inequality \eqref{M_eta < c M} holds true.  \\
2. We have
$$M_{\lambda,\eta}(v) - M(v) = \bigg(\frac{M_{\lambda,\eta}(v)}{p_2^{\lambda,\eta}(v)}-\frac{M(v)}{p_2^0(v)}\bigg)p_2^{\lambda,\eta}(v)+\frac{M(v)}{p_2^0(v)}\big(p_2^{\lambda,\eta}(v)-p_2^0(v)\big) .$$
So, for $|v|\leqslant s_0\eta^{-\frac{1}{3}}$, since $p_2^{\lambda,\eta}(v)=p_2^0(v)$ then,
$$|M_{\lambda,\eta}(v) - M(v)| \leqslant p_2^0(v)\left\|\frac{M_{\lambda,\eta}}{p_2^{\lambda,\eta}}-\frac{M}{p_2^0}\right\|_{\infty} \leqslant \left\|\frac{M_{\lambda,\eta}}{p_2^{\lambda,\eta}}-\frac{M}{p_2^0}\right\|_{\infty} \underset{\eta \rightarrow 0}{\longrightarrow}0.$$
For $|v|\geqslant s_0\eta^{-\frac{1}{3}}$, we have $p_2^{\lambda,\eta}(v)\leqslant p_2^0(v)\leqslant\eta^{\frac{\gamma}{3}}$, then
$$|M_{\lambda,\eta}(v) - M(v)| \leqslant \eta^{\frac{\gamma}{3}}\bigg(2+\left\|\frac{M_{\lambda,\eta}}{p_2^{\lambda,\eta}}-\frac{M}{p_2^0}\right\|_{\infty}\bigg) \underset{\eta \rightarrow 0}{\longrightarrow}0.$$
Then \eqref{cv dans L^2} is obtained by Lebesgue's theorem.  \\
3. Let $\nu =\frac{1}{2}(\gamma-\frac{1}{2})$. We have  $2\gamma-\nu-1>0$. Then, $\langle v \rangle^{\nu-2\gamma} \in L^1(\mathbb{R})$ and $|v|\langle v \rangle^{-\nu}\leqslant \eta^{\frac{\nu-1}{3}}$ for $|v|\leqslant\eta^{-\frac{1}{3}}$. Now, since $|M_{\lambda,\eta}(v)|\lesssim p^{\lambda,\eta}_2(v)$ with $p^{\lambda,\eta}_2(v)=\langle v \rangle^{-\gamma}$ for $|v|\leqslant s_0\eta^{-\frac{1}{3}}$,  then
$$\bigg|\int_\mathbb{R}\eta^{\frac{1}{3}}v M_{\lambda,\eta}(v)M(v)\mathrm{d}v\bigg|\lesssim \int_{|v|\leqslant s_0\eta^{-\frac{1}{3}}}{\eta^{\frac{1}{3}}}|v|\langle v \rangle^{-\nu+\nu-\gamma}M(v)\mathrm{d}v+\int_{|v|\geqslant s_0\eta^{-\frac{1}{3}}} {\eta^{\frac{1}{3}}}|v|p^{\lambda,\eta}_2(v)M(v)\mathrm{d}v. $$
Therefore, after making the change of variable $v=\eta^{-\frac{1}{3}}s$ in the second integral, we get
\begin{align*}
\bigg|\int_\mathbb{R}\eta^{\frac{1}{3}}v M_{\lambda,\eta}(v)M(v)\mathrm{d}v\bigg|\lesssim \bigg[\eta^{\frac{\nu}{3}}\int_{\mathbb{R}}\langle v \rangle^{\nu-2\gamma}\mathrm{d}v +\eta^{\frac{2\gamma-1}{3}}\int_{|s|\geqslant s_0} |s|^{1-\gamma}|\mathfrak{a}_{\lambda}(|s|)|\mathrm{d}s\bigg]\underset{\eta \rightarrow 0}{\longrightarrow 0} ,
\end{align*}
with $\mathfrak{a}_{\lambda}(s)=\mathfrak{a}_{\lambda,\eta}(v\eta^{-\frac{1}{3}})$ does not depend on $\eta$ after rescaling.
\ep
\begin{corollary}[Rescaled solution]\label{Propr de H_lambda,eta}
We define the function $H_{\lambda,\eta}$ for all $s$ in $\mathbb{R}$ by 
\begin{equation}\label{H_lambda,eta}
H_{\lambda,\eta}(s) := \eta^{-\frac{\gamma}{3}}M_{\lambda,\eta}(\eta^{-\frac{1}{3}}s) .
\end{equation}
Then, $H_{\lambda,\eta}$ satisfies the rescaled equation
$$\big[-\partial^2_s+ \frac{\gamma(\gamma+1)}{|s|_\eta^2} + \mathrm{i} s -\lambda \big]H_{\lambda,\eta}(s)= -V_\eta(s)H_{\lambda,\eta}(s)-\eta^{-\frac{2+\gamma}{3}}b(\lambda,\eta)\Phi(\eta^{-\frac {1}{3}}s),\quad  s\in \RR  .$$
Moreover, the following estimates hold
\begin{enumerate}
\item For all $|s|\leqslant s_0$ 
\begin{equation}\label{H_lambda,eta < s^-gamma}
|H_{\lambda,\eta}(s)| \lesssim |s|_\eta^{-\gamma} \leqslant |s|^{-\gamma}  .
\end{equation}
\item For all $|s|\geqslant s_0$
\begin{equation}\label{H_lambda,eta < a_lambda}
|H_{\lambda,\eta}(s)| \lesssim \left\{\begin{array}{l} |\mathfrak{a}_{\lambda}(s)| ,   \quad s \geqslant s_0 , \\
|\mathfrak{b}_{\lambda}(s)| , \quad s \leqslant -s_0 ,
\end{array}\right.
\end{equation}
\end{enumerate}
where $|s|_\eta := (\eta^\frac{2}{3}+s^2)^\frac{1}{2}$ and $V_\eta(s) := \big[\gamma(\gamma+1)s^2-\gamma\eta^\frac{2}{3}\big]|s|_\eta^{-4}$.
\end{corollary}
\bp The proof of this Corollary is a direct consequence of the previous one, by making the change of variable $v=\eta^{-\frac{1}{3}}s$.
\ep
\subsection{Existence of an eigen couple $(\mu(\eta),M_{\mu,\eta})$ for the complete operator}
The purpose of this section is, first, to find a solution for equation \eqref{M_mu,eta}, which amounts to showing the existence of a $\lambda$,  function of $\eta$, such that the additional term $\langle M_{\lambda,\eta}-M,\Phi\rangle=0$. We prove it again using the implicit function theorem. In a second step, we will compute the eigenvalue $\mu =\lambda\eta^{\frac{2}{3}}$, and for this, we will establish some estimates on the solutions of the equations of $M_{\lambda,\eta}$ and $H_{\lambda,\eta}$ respectively.
We summarize these two results in the following two propositions:
\begin{proposition}[Constraint]\label{contrainte}
Define $$ B(\lambda,\eta):= \eta^{-\frac{2}{3}} b(\lambda,\eta).$$
\begin{enumerate}
\item The expression of $B(\lambda,\eta)$ is given by
\begin{equation}\label{b(lambda,eta)}
B(\lambda,\eta)= \eta^{-\frac{2}{3}}\langle N_{\lambda,\eta},\Phi\rangle = \int_\mathbb{R}(\lambda-\mathrm{i}\eta^{\frac{1}{3}} v)M_{\lambda,\eta}(v)M(v)\mathrm{d}v .
\end{equation}
\item The $\eta$ order of the coefficient in front of $\lambda$ in the expansion on $\lambda$ of $B(\lambda,\eta)$ is given by
\begin{equation}\label{B(lambda,0)}
\underset{\eta \rightarrow 0}{\lim} B(\lambda,\eta) = \lambda\int_\mathbb{R} M^2(v)\mathrm{d}v .
\end{equation}
\item There exists $\tilde\eta_0, \tilde\lambda_0>0$ small enough, a function $\tilde{\mu}: \{|\eta|\leqslant\tilde\eta_0\} \longrightarrow \{|\lambda|\leqslant\tilde\lambda_0\}$ such that:\\
for all $(\lambda,\eta)\in \{|\eta|<\tilde\eta_0\}\times\{|\lambda|<\tilde\lambda_0\}$,  $\lambda = \tilde{\mu}(\eta)$ and the constraint is satisfied: $$B(\lambda,\eta)= B(\tilde{\mu}(\eta),\eta)=0.$$
\end{enumerate}
Therefore, $\mu(\eta) = \eta^{\frac{2}{3}} \tilde \mu(\eta)$ is the eigenvalue associated to the eigenfunction $M_\eta: = M_{\tilde \mu(\eta),\eta}$ for the operator $\mathcal{L}_\eta$, and the couple $\big(\mu(\eta),M_\eta \big)$ is a solution to the spectral problem \eqref{M_mu,eta}. 
\end{proposition}
\begin{proposition}[Approximation of the eigenvalue]\label{vp} Let $\alpha:=\frac{2\gamma+1}{3}$ for all $\gamma \in ]\frac{1}{2},\frac{5}{2}[$.
The eigenvalue $\mu(\eta)$ satisfies
\begin{equation}\label{mu(eta)}
\mu(\eta)=\overline \mu (-\eta)= \kappa |\eta|^{\alpha}\big(1+O(|\eta|^\alpha)\big) ,
\end{equation}
where  $\kappa$ is a positive constant given by
\begin{equation}\label{kappa}
\kappa := - 2 C^2_\beta \int_{0}^\infty s^{1-\gamma} \Imm H_0(s) \mathrm{d}s ,
\end{equation}
and where $H_0$ is the unique solution to
\begin{equation}\label{eq de H_0}
\big[-\partial^2_s+ \mathrm{i} s+\frac{\gamma(\gamma+1)}{s^2}\big]H_0(s)=0, \ s\in \mathbb{R}^* ,
\end{equation}
satisfying \begin{equation}\label{condition H_0}
\int_{|s|\geqslant 1}|H_0(s)|^2\mathrm{d}s <+\infty \  \mbox{ and } \  H_0(s)\underset{0}{\sim} |s|^{-\gamma} .
\end{equation}
\end{proposition}

\noindent \bpp \ref{contrainte}.\\
1. The first point is obtained by multiplying the equation \eqref{eq de M_eta + cont.} by $M$, integrating it twice by part and using the fact that $[-\partial_v^2+W(v)]M=0 $ and $\langle \Phi,M\rangle =1$.  \\
2. We obtain the limit \eqref{B(lambda,0)} by the last two points of the Corollary \ref{Propr de M_lambda,eta}.  \\
3. The proof of this point is an immediate consequence of the IFT applied to the function $B$ around the point $(0,0)$.
\epp

\noindent In order to get the Proposition \ref{vp}, we need to prove the following two lemmas:
The first one gives estimates on $M_{0,\eta}$ and $H_{0,\eta}$.
\begin{lemma}\label{lemme N_0,eta}
For all $\gamma> 1$ ons has
\begin{equation}\label{estimation de N_0,eta v petit}
|M_{0,\eta}(v)-M(v)|\lesssim \eta ,\quad \forall v \in [-v_0,v_0] .
\end{equation}
Moreover, for large velocities
\begin{equation}\label{estimation de N_0,eta v grand}
|M_{0,\eta}(v)-M(v)|\lesssim \eta \langle v\rangle^{3-\gamma},\ \forall |v| \in [v_0,s_0\eta^{-\frac{1}{3}}] .
\end{equation}
Therefore,
\begin{equation}\label{H_0,eta-s^-gamma,eta<s^(3-delta_2-gamma)}
\big|H_{0,\eta}(s)-|s|_\eta^{-\gamma}\big| \lesssim \ |s|_\eta^{3-\gamma}  \ \leqslant \langle s\rangle^{3-\gamma} ,  \ \ \forall |s| \in [0,s_0] .
\end{equation}
\end{lemma}
The second one gives the formula of the diffusion coefficient.
\begin{lemma}\label{lim kappa_eta}
\item \begin{enumerate}
\item The small velocities don't participate to the limit in the approximation of $\mu(\eta)$:
\begin{equation}\label{petites vitesses}
\underset{\eta \rightarrow 0}{\lim} \eta^{-\frac{2(\gamma-1)}{3}} \int_{|v|\leqslant v_0} v M_{0,\eta}(v) M(v) \mathrm{d}v = 0 .
\end{equation}
\item For large velocities, we have the following limit:
\begin{equation}\label{lim kappa_eta = kappa}
\underset{\eta \rightarrow 0}{\lim}\ \mathrm{i} \eta^{-\frac{2(\gamma-1)}{3}} \int_{|v| \geqslant v_0} v M_{0,\eta}(v) M(v) \mathrm{d}v = - 2 \int_{0}^\infty s^{1-\gamma} \Imm H_0(s) \mathrm{d}s , 
\end{equation}
where $H_0$ is the unique solution to \eqref{eq de H_0} satisfying the conditions \eqref{condition H_0}.
\end{enumerate}
\end{lemma}
\bpl \ref{lemme N_0,eta}.  Recall that $N_{0,\eta}:=M_{0,\eta}-M$ and it satisfies the equation 
\begin{equation}\label{N_0,eta}
[-\partial^2_v+W(v)]N_{0,\eta}(v)= - \mathrm{i} \eta v [N_{0,\eta}(v)+M(v)] -\langle N_{0,\eta},\Phi\rangle\Phi(v), \ \forall v \in \RR .
\end{equation} 
Thanks to the symmetry $\bar N_{0,\eta}(-v)=N_{0,\eta}(v)$, we establish the inequalities \eqref{estimation de N_0,eta v petit} and \eqref{estimation de N_0,eta v grand} on $[0,s_0\eta^{-\frac{1}{3}}]$. By writing the solution of the equation \eqref{N_0,eta} in the basis of solutions of $[-\pa_v^2+W(v)]f=0$, which is given by $\{M,Z\}$, we get:
\begin{align*}
N_{0,\eta}(v)=&\bigg(c_1-\int_0^v [\mathrm{i} \eta wZ (N_{0,\eta}+M)+b(0,\eta)\Phi(w)]\mathrm{d}w\bigg)M(v) \\
&+\bigg(c_2 + \int_0^v [\mathrm{i} \eta wM(N_{0,\eta}+M)+b(0,\eta)\Phi(w)]\mathrm{d}w\bigg)Z(v) ,
\end{align*}
where $c_1$ and $c_2$ are two complex constants to determine and $b(0,\eta):=\langle N_{0,\eta},\Phi\rangle$.
Now, by using the condition $N_{0,\eta}(0)=0$, we get $c_1=0$, and since $|N_{0,\eta}(\eta^{-\frac{1}{3}})| \lesssim M(\eta^{-\frac{1}{3}}) \sim \eta^{\frac{\gamma}{3}}$ (that we get from \eqref{M_eta < c M} for $v=\eta^{-\frac{1}{3}}$) then, 
$$c_2= -\int_0^\infty [\mathrm{i} \eta wM(N_{0,\eta}+M)+b(0,\eta)\Phi(w)]\mathrm{d}w  , $$ 
otherwise we will have $|N_{0,\eta}(\eta^{-\frac{1}{3}})| \sim \eta^{-\frac{\gamma-2}{3}}$, which contradicts the fact that $|N_{0,\eta}(\eta^{-\frac{1}{3}})| \lesssim \eta^{\frac{\gamma}{3}}$ for all $\gamma \in (\frac{1}{2},\frac{5}{2})$. Therefore,
\begin{align*}
N_{0,\eta}(v)=&-\mathrm{i} \eta\bigg[M(v)\int_0^v wZ (N_{0,\eta}+M)\mathrm{d}w+Z(v)\int_v^\infty wM(N_{0,\eta}+M)\mathrm{d}w\bigg] \\
&- b(0,\eta)\bigg[M(v)\int_0^v \Phi Z \mathrm{d}w + Z(v) \int_v^\infty \Phi M \mathrm{d}w\bigg] .
\end{align*}
By \eqref{b(lambda,eta)}, $b(0,\eta) = -\mathrm{i}\eta \int_\RR wM_{ 0,\eta}M\mathrm{d}w $, and since $|N_{0,\eta}(v)| \lesssim M(v)$ then, by dividing the last equality by $\eta \langle v \rangle^{3-\gamma}$, we get
$$
\bigg|\frac{N_{0,\eta}(v)}{\eta \langle v \rangle^{3-\gamma}}\bigg| \lesssim \langle v \rangle^{-3}\int_0^v [wM+\Phi]Z \mathrm{d}w+\langle v \rangle^{\gamma-3}Z(v)\int_v^\infty [wM+\Phi]\mathrm{d}w .\\
$$
We have $Z(v)=M(v)\int_0^v\frac{\mathrm{d}w}{M^2(w)}\lesssim \langle v \rangle^{\gamma+1}$ and $\Phi(v) \leqslant p^0_1(v) \leqslant \langle v \rangle^{-\gamma-2}$. Thus, $\big|\frac{N_{0,\eta}(v)}{\eta \langle v \rangle^{3-\gamma}}\big| \lesssim 1$ for all $v \in \RR^+$, in particular for $v \in [0,v_0]$ we get \eqref{estimation de N_0,eta v petit} and for $v \in [0,s_0 \eta^{-\frac{1}{3}}]$ we get \eqref{estimation de N_0,eta v grand}. \epl  \\

\noindent \bpl \ref{lim kappa_eta}. \\
First of all, since $\bar{M}_{0,\eta}(-v)=M_{0,\eta}(v)$ and $M(-v)=M(v)$ for all $v\in\RR$, then
$$ - \mathrm{i} \int_\mathbb{R} vM_{0,\eta}(v)M(v)\mathrm{d}v = 2 \int_0^\infty v \Imm M_{0,\eta}(v)M(v)\mathrm{d}v = 2 \int_0^\infty v \Imm \big(M_{0,\eta}(v)-M(v)\big)M(v)\mathrm{d}v .$$ 
1. Recall that $\alpha:=\frac{2\gamma+1}{3}$. Then,  $1-\alpha = \frac{2(\gamma-1)}{3}$.  We have thanks to \eqref{M_eta < c M} and \eqref{estimation de N_0,eta v petit}:
$$
\eta^{1-\alpha}\bigg|\int_{|v|\leqslant v_0} vM_{0,\eta}(v)M(v)\mathrm{d}v\bigg| \lesssim \left\{\begin{array}{l} \ds \eta^{\frac{2(1-\gamma)}{3}}\int_{|v|\leqslant v_0} \langle v\rangle^{1-2\gamma}\mathrm{d}v ,  \quad \gamma\in (\frac{1}{2},1) ,\\
 \\
 \ds  \eta^{\frac{5-2\gamma}{3}}\int_{|v|\leqslant v_0} \langle v\rangle^{4-\gamma} \mathrm{d}v ,  \qquad \gamma\in (1,\frac{5}{2}) .
\end{array}\right. 
$$
Hence  \eqref{petites vitesses} holds true for all $\gamma \in (\frac{1}{2},\frac{5}{2})$, since $2-\alpha=(5-2\gamma)/3 >0$ for $\gamma>1$.  The case $\gamma=1$ is done by Lebesgue's theorem thanks to \eqref{limite simple M_eta(v) = M(v)} and the domination of $v|M_{0,\eta}|M$ by $v\langle v\rangle^{-2\gamma} \in L^1(0,v_0)$ thanks to \eqref{M_eta < c M}.\\
2. For the the integral on $[v_0,+\infty)$, we split it into two parts as follows:
$$
2 \eta^{1-\alpha} \bigg[\int_{v_0}^{s_0\eta^{-\frac{1}{3}}} v\Imm\big(M_{0,\eta}(v)-M(v)\big)M(v)\mathrm{d}v+
\int_{s_0\eta^{-\frac{1}{3}}}^\infty v\Imm M_{0,\eta}(v) M(v)\mathrm{d}v\bigg] .
$$
In order to compute the limit of the previous expressions, we proceed to a change of variable $v=\eta^{-\frac{1}{3}}s$, which means that we need to compute
$$ 
\underset{\eta \rightarrow 0}{\lim} \bigg[\int_{\eta^{\frac{1}{3}}v_0 }^{s_0} s|s|_\eta^{-\gamma}\Imm\big[H_{0,\eta}(s)-|s|_\eta^{-\gamma}\big] \mathrm{d}s+\int_{s_0}^\infty s|s|_\eta^{-\gamma} \Imm H_{0,\eta}(s)\mathrm{d}s\bigg]  .
$$
For that purpose, we will use the ``weak-strong" convergence in the Hilbert space $L^2(0,\infty)$.  First, recall the following estimates given by \eqref{H_lambda,eta < s^-gamma}, \eqref{H_lambda,eta < a_lambda} and \eqref{H_0,eta-s^-gamma,eta<s^(3-delta_2-gamma)}:  \\
$\bullet$ For $s \in (0,s_0)$,
$$
\big| H_{0,\eta}(s)-|s|_\eta^{-\gamma}\big| \lesssim \left\{\begin{array}{l} |s|_\eta^{-\gamma} ,  \quad \gamma\in (\frac{1}{2},1] \\
 \\
 |s|_\eta^{3-\gamma} ,  \ \ \gamma\in (1,\frac{5}{2})
\end{array}\right. \lesssim \left\{\begin{array}{l} |s|^{-\gamma} , \ \ \quad \gamma\in (\frac{1}{2},1] ,\\
 \\
 \langle s\rangle^{3-\gamma} ,  \quad \gamma\in (1,\frac{5}{2}) .
\end{array}\right.
$$
$\bullet$ For $s \geqslant s_0$ and for all $\gamma \in (\frac{1}{2},\frac{5}{2})$,
$$\big| H_{0,\eta}(s)\big| \lesssim \big| \mathfrak{a}_0(s)\big| .$$
Hence, the two sequences $(H_{0,\eta})_\eta$ and $(\mathcal{H}_\eta)_\eta$ are uniformly bounded in $L^1(0,\infty)$ and $L^2(0,\infty)$ respectively, where $\mathcal{H}_\eta$ is defined by
$$
\mathcal{H}_\eta(s):=  \left\{\begin{array}{l} s^{\frac{1}{2}} \Imm H_{0,\eta}(s)=s^{\frac{1}{2}} \Imm\big[H_{0,\eta}(s)-|s|_\eta^{-\gamma}\big],   \qquad \gamma\in (\frac{1}{2},1], 0<s \leqslant s_0 , \\
 \\
s^{-1} \Imm H_{0,\eta}(s)= s^{-1} \Imm \big[H_{0,\eta}(s)-|s|_\eta^{-\gamma}\big] , \quad \ \ \gamma\in (1,\frac{5}{2}), 0<s \leqslant s_0 , \\
 \\
 s \Imm H_{0,\eta}(s)   \qquad  \qquad \mbox{  for all } \gamma\in (\frac{1}{2},\frac{5}{2}) \mbox{ and } s \geqslant s_0 .
\end{array}\right.
$$
Now, since the sequence $\mathcal{H}_\eta$ is bounded in $L^2(0,\infty)$, uniformly with respect to $\eta$ then, up to a subsequence, that we will again denote by $\mathcal {H}_\eta$, we have $\mathcal{H}_\eta$ converges weakly in $L^2(0,\infty)$, thus converges in $\mathcal{D}'(0,\infty)$. Let's identify this limit that we denote by $\mathcal{H}_0$.
We have on the one hand, $H_{0,\eta}$ converges to $H_0$ in $\mathcal{D}'(0,\infty)$. Indeed,  recall that $H_{0,\eta}$ satisfies the equation
$$
\big[-\partial^2_s+\frac{\gamma(\gamma+1)}{|s|^2_\eta} + \mathrm{i} s \big] H_{0,\eta}(s)= \eta^{\frac{2}{3}}\frac{\gamma(\gamma+2)}{|s|^4_\eta} H_{0,\eta}(s) - \eta^{-\frac{2+\gamma}{3}} b(0,\eta) \Phi(\eta^{-\frac{1}{3}}s) .
$$
Let $\varphi \in \mathcal{D}(0,\infty)$. Then,  by multiplying the previous equation by $\varphi$ and by integrating it by parts, we obtain
$$ 
\int_0^\infty \bigg[-\partial_s^2+ \frac{\gamma(\gamma+1)}{|s|^2_\eta}+\mathrm{i} s \bigg]\varphi(s) H_{0,\eta}(s)\mathrm{d}s = \int_0^\infty \bigg[\eta^{\frac{2}{3}}\frac{\gamma(\gamma+2)}{|s|^4_\eta}H_{0,\eta}(s) - \eta^{-\frac{\gamma}{3}} B(0,\eta)\Phi(\eta^{-\frac{1}{3}}s)\bigg]\varphi(s)\mathrm{d}s .
$$
By \eqref{b(lambda,eta)} and \eqref{int eta^1:3 v M_eta M->0}, $B(0,\eta) \to 0$ when $\eta \to 0$ and since $\Phi(\eta^{-\frac{1}{3}}s) \lesssim p^{0,\eta}_1(\eta^{-\frac{1}{3}}s) \leqslant \eta^{\frac {2+\gamma+\delta}{3}}|s|^{-2-\gamma-\delta}$ then, by Lebesgue's theorem $H_{0,\eta}$ converges to $H_0$, in  $\mathcal{D}'(0,\infty)$, the unique solution to 
$$
\big[-\partial_s^2+  \frac{\gamma(\gamma+1)}{s^2}+ \mathrm{i} s \big]H_0(s)= 0 .
$$
The uniqueness of $H_0$ comes from the fact that the previous equation admits a unique solution in $L^2(1,\infty)$, up to a multiplicative constant \cite{LebPu} and that $H_0(s) \underset{0}{\sim} s^{-\gamma}$ which is obtained by passing to the limit, in $\eta$, in the inequality $| H_{0,\eta}(s_1) / |s_1|^{-\gamma}_\eta - 1 | \leqslant \| M_{0,\eta} / p^ {0,\eta}_2-M / p^0_2 \|_\infty$ for $s_1>0$ small enough. \\
On the other hand, since $s^{-\frac{1}{2}}$ and $s$ are in $C^\infty(0,\infty)$ and $H_{0,\eta} \to H_0$ in $\mathcal{D}'(0,\infty)$ with $H_0$ unique.  Then, $\mathcal{H}_\eta$ converges to $\mathcal{H}_0$ in $\mathcal{D}'(0,\infty)$ and therefore $\mathcal{H}_\eta \rightharpoonup \mathcal{H}_0$ weakly in $L^2(0,\infty)$ with $\mathcal{H}_0$ unique and given by 
$$
\mathcal{H}_0(s):=   \left\{\begin{array}{l} s^{\frac{1}{2}} \Imm H_0(s),   \qquad \gamma\in (\frac{1}{2},1], 0<s \leqslant s_0 ,\\
 \\
 s^{-1} \Imm H_0(s) , \quad \ \ \gamma\in (1,\frac{5}{2}), 0<s \leqslant s_0 , \\
 \\
 s \Imm H_0(s)  \qquad   \qquad \mbox{  for all } \gamma\in (\frac{1}{2},\frac{5}{2}) \mbox{ and } s \geqslant s_0 .
\end{array}\right.
$$
Moreover, thanks to the uniqueness of this limit,  the whole sequence converges.
Finally, we conclude by passing to the limit in the scalar product $\langle \mathcal{H}_\eta, \mathcal{I}_\eta\rangle$, where $\mathcal{I}_\eta$ definded by 
$$
\mathcal{I}_\eta:=  \left\{\begin{array}{l} s^{\frac{1}{2}} |s|_\eta^{-\gamma},  \qquad \gamma\in (\frac{1}{2},1], 0<s \leqslant s_0 , \\
 \\
 s^{2} |s|_\eta^{-\gamma} , \quad \quad \gamma\in (1,\frac{5}{2}), 0<s \leqslant s_0 ,\\
 \\
|s|_\eta^{-\gamma} ,   \qquad \quad  \gamma\in (\frac{1}{2},\frac{5}{2}) , s \geqslant s_0 ,
\end{array}\right.  
$$
converges strongly in $L^2(0,\infty)$ to 
$$ 
 \mathcal{I}_0 =  \left\{\begin{array}{l} s^{\frac{1}{2}-\gamma},  \qquad \gamma\in (\frac{1}{2},1], 0<s \leqslant s_0 ,\\
 \\
  s^{2-\gamma}  , \qquad \gamma\in (1,\frac{5}{2}), 0<s \leqslant s_0 , \\
 \\
  s^{-\gamma}  ,  \quad  \quad \ \gamma\in (\frac{1}{2},\frac{5}{2}) ,  s \geqslant s_0 .
\end{array}\right.
$$
Hence the limit \eqref{lim kappa_eta = kappa} holds true.
\epl

\noindent \bpp \ref{vp}.  By doing an expansion in $\lambda$ for $B$ and by Proposition \ref{contrainte}, we get
$$
B(\lambda,\eta)=\eta^{-\frac{2}{3}}b(\lambda, \eta)=\eta^{-\frac{2}{3}}b(0,\eta)+ \lambda \int_\RR M_{0,\eta} M \mathrm{d}v+ O(\lambda^2) .
$$ 
Then, for $\lambda=\tilde{\mu}(\eta)$ and since $B(\tilde{\mu}(\eta),\eta)=0$, we obtain
$$ 
\tilde{\mu}(\eta)= - \eta^{-\frac{2}{3}}b(0,\eta)\bigg(\int_\RR M_{0,\eta} M \mathrm{d}v\bigg)^{-1}+o\big(\eta^{-\alpha}b(0,\eta)\big) ,
$$ 
which implies that 
$$\eta^{-\alpha}\mu(\eta)=\eta^{\frac{2}{3}-\alpha}\tilde{\mu}(\eta)= -\eta^{-\alpha}b(0,\eta)\bigg(\int_\RR M_{0,\eta} M \mathrm{d}v\bigg)^{-1} .
$$
By \eqref{cv dans L^2} and \eqref{lim kappa_eta = kappa}, $\underset{\eta \rightarrow 0}{\lim} \int_{\mathbb{R}} M_{0,\eta}(v)M(v)\mathrm{d}v = \|M\|^2_2$ and $\underset{\eta \rightarrow 0}{\lim} \eta^{-\alpha}b(0,\eta) = 2\int_0^\infty s^{1-\gamma}\Imm H_0(s)\mathrm{d}s$ respectively. Hence,  $\underset{\eta \rightarrow 0}{\lim} \eta^{-\alpha}\mu(\eta) = \kappa$.
For $\eta\in [-\eta_0, 0]$, the symmetry $ \mu(\eta) = \overline \mu (-\eta) $ holds by complex conjugation on the equation. Thus, the proof of Proposition \ref{vp} is complete.
\epp
\begin{remark}
We do the same calculations as in \cite{LebPu} to compute the coefficient $\kappa$. Also, since it is given by the same integral formula \eqref{kappa} then, we have $\kappa > 0$.
\end{remark}

\noindent \textbf{Proof of the main Theorem \ref{main}.} The existence and uniqueness of the eigen-solution $\big(\mu(\eta),M_\eta\big)$ is given by  Proposition \ref{contrainte}. The first point, \eqref{item 1 thm}, is given by the Corrolary \ref{Propr de M_lambda,eta} and finally, the second point is given by the Proposition \ref{vp}.  \hfill $\square$

\end{document}